\documentclass[11pt]{article}
\newcommand{\documentdate}{15 I 2025}

\usepackage{a4wide,latexsym,amsmath,varioref,graphicx,xcolor,amssymb,dsfont}

\title{A Stochastic Objective-Function-Free Adaptive Regularization
  Method with Optimal Complexity} 

\author{	
   S. Gratton%
   \thanks{Universit\'e de Toulouse, INP, IRIT, Toulouse, France. Email:
     serge.gratton@enseeiht.fr. Work partially supported by 3IA Artificial and
     Natural Intelligence Toulouse Institute (ANITI), French "Investing for the Future
     - PIA3" program under the Grant agreement ANR-19-PI3A-0004"}, 
   ~S. Jerad%
   \thanks{Universit\'e de Toulouse, INP, IRIT, Toulouse,
     France. S. Jerad now works in the Mathematical Institute,
     University of Oxford, United Kingdom. Email: sadok.jerad@maths.ox.ac.uk}
   ~and Ph. L. Toint%
   \thanks{NAXYS, University of Namur, Namur, Belgium. Email:
     philippe.toint@unamur.be. Partly supported by ANITI.}
}

\newcommand{\beqn}[1]{\begin{equation}\label{#1}}
\newcommand{\eeqn}{\end{equation}}
\newcommand{\req}[1]{(\ref{#1})}
\newcommand{\ms}{\;\;\;\;}
\newcommand{\tim}[1]{\;\; \mbox{#1} \;\;}

\newtheorem{theorem}{Theorem}[section]
\newtheorem{lemma}[theorem]{Lemma}
\newtheorem{corollary}{Corollary}[section]

\newcommand{\numsection}[1]{\section{#1}\setcounter{equation}{0}}
\renewcommand{\theequation}{\arabic{section}.\arabic{equation}}

\newcounter{algo}[section]
\renewcommand{\thealgo}{\thesection.\arabic{algo}}
\newcommand{\llem}[2]{\vspace{\baselineskip} 
\noindent\framebox[\textwidth]{\parbox{0.95\textwidth}{
\begin{lemma} \label{#1} \rm #2 \end{lemma} } } \vspace{\baselineskip} }
\newcommand{\algo}[3]{\refstepcounter{algo}
\begin{center}\begin{figure}[htbp]
\framebox[\textwidth]{
\parbox{0.95\textwidth} {\vspace{\topsep}
{\bf Algorithm \thealgo : #2}\label{#1}\\
\vspace*{-\topsep} \mbox{ }\\
{#3} \vspace{\topsep} }}
\end{figure}\end{center}}
\newcommand{\bpr}{{\bf Proof.} \hspace{1.5mm}}
\newcommand{\epr}{\hfill $\Box$ \vspace*{1em}}
\newcommand{\proof}[1]{
\begin{list}{}{
\setlength{\topsep}{0.0pt}
\setlength{\partopsep}{0.0pt}
\setlength{\leftmargin}{0.025\textwidth}
\setlength{\rightmargin}{0.5\leftmargin}
\setlength{\labelwidth}{0.5\leftmargin}
\setlength{\labelsep}{0.25\leftmargin}}
\item \bpr #1 \epr \noindent
\end{list}}
\newcommand{\lthm}[2]{\vspace{\baselineskip} 
\noindent\framebox[\textwidth]{\parbox{0.95\textwidth}{
\begin{theorem} \label{#1} \rm #2 \end{theorem} } } \vspace{\baselineskip} }
\newcommand{\lcor}[2]{\vspace{\baselineskip} 
\noindent\framebox[\textwidth]{\parbox{0.95\textwidth}{
\begin{corollary} \label{#1} \rm #2 \end{corollary} } } \vspace{\baselineskip}
}
\newcommand{\ii}[1]{\{ 1, \ldots, #1 \}}
\newcommand{\iiz}[1]{\{ 0, \ldots, #1 \}}
\newcommand{\iibe}[2]{\{ #1, \ldots, #2 \}}

\newcommand{\E}[1]{\mathbb{E}\left[#1\right]}
\newcommand{\Ek}[2]{\mathbb{E}_{#2} \left[#1\right]}
 
\newcommand{\calB}{{\cal B}} 
\newcommand{\calF}{{\cal F}} 
\newcommand{\calO}{{\cal O}}

\newcommand{\twoppp}{2^\sfrac{1}{p}}

\renewcommand{\Re}{\hbox{I\hskip -2pt R}}

\newcommand{\sfrac}[2]{{\scriptstyle \frac{#1}{#2}}}
\newcommand{\half}{\sfrac{1}{2}}
\newcommand{\eqdef}{\stackrel{\rm def}{=}}
\newcommand{\bigsum}{\displaystyle \sum}
\newcommand{\kap}[1]{\kappa_{\mbox{\tiny #1}}}
\newcommand{\khigh}{\kappa_{\mbox{\tiny high}}}
\newcommand{\al}[1]{{\footnotesize{\sf #1}}}
\newcommand{\tal}[1]{{\normalsize {\sf #1}}}

\newcommand{\phil}[1]{{\color{teal}#1}}

\newcommand{\overder}[2]{\overline{\nabla_x^{#1} f}(X_{#2})}
\newcommand{\indica}[1]{\mathds{1}_{#1}}

\newcommand{\pplusonefacpow}{(p+1)!^\sfrac{p+1}{p}}

\newcommand{\sumskzeroterm}{  \twoppp \left(\frac{\kappa_D}{\sigma_0^\sfrac{p+1}{p}} + \frac{\E{\|G_0\|^{\sfrac{p+1}{p}}}}{\sigma_{0}^\sfrac{p+1}{p}} \right)}
\newcommand{\bigNred}{\textcolor{red}{N}}
\newcommand{\comment}[1]{}

\date{\documentdate}

\newcounter{assumption}

\renewcommand{\theassumption}{AS.\arabic{assumption}}

\newenvironment{assumption}[1]
{
    \refstepcounter{assumption}
    \noindent \textbf{\theassumption} \label{#1}
}
{}

\begin{document}
	
\maketitle

\begin{abstract}
A fully stochastic $p$th-order adaptive-regularization method for
unconstrained nonconvex optimization is presented which never computes
the objective-function value, but yet achieves the optimal
$\mathcal{O}(\epsilon^{-(p+1)/p})$ complexity bound for finding
first-order critical points. When stochastic gradients and Hessians
are considered, we recover the optimal
$\mathcal{O}\left(\epsilon^{-3/2}\right)$ bound for finding
first-order critical points.  The method is noise-tolerant and the
inexactness conditions required for convergence depend on the history
of past steps. Applications to cases where derivative evaluation is
inexact and to minimization of finite sums by sampling are discussed.
Numerical experiments on large binary classification problems
illustrate the potential of the new method.
\end{abstract}

{\small
\textbf{Keywords:} stochastic optimization, adaptive regularization methods,
evaluation complexity, \\Objective-Function-Free-Optimization (OFFO), nonconvex optimization.
}
\numsection{Introduction}

Adaptive gradient methods such as Adam \cite{kingma2014adam}, Adagrad
\cite{duchi11adagrad}, or AMSGrad \cite{ReddiKK18} have become the
workhorse of large-scale stochastic optimization, especially when
training artificial neural networks. Given their remarkable empirical
results, various complexity analyses have been developed in the
stochastic regime
\cite{LiOrabona2018,defossez2022a,ZhouChenTangYangCaoGu20,Fawetal22}
or the deterministic one \cite{WardWuBott19,GraJerToin24}. A notable
feature of these methods is that they do not compute the function
value or an approximation thereof, making them part of OFFO (Objective
Free Function Optimization) methods
\cite{WardWuBott19,GraJerToin24}. Examples of such methods are 
adaptive gradient methods, which use only the current and past
gradients, and are ubiquitous in machine learning due to their
performance on \phil{noisy} large-dimensional problems. The
  motivation for these linesearch-based methods stems from the theoretical
  \cite{Berahas2021, Cartis2017, Paquette2020}
  and practical
  observation that the accuracy needed for a noisy function value is
  significantly higher than that required for the gradient. When
  function and gradient values are obtained by sampling, the size of
  the necessary sample for the function value is typically much larger
  than for the gradient, making the use of stochastic subsamples
  for the estimation of function values both numerically intensive and
  impractical. Adaptive gradient algorithms however remain
  theoretically limited by a complexity bound in
$\mathcal{O}\left(\epsilon^{-2}\right)$ for
finding an $\epsilon$-approximate first order solution
\cite{WardWuBott19,GraJerToin24}.
Moreover, even higher-order schemes such as trust region
\cite{Blanchet2019,Chen2017}, or adaptive regularization methods
\cite{Bellavia2019,Bellavia2022,toint2020stochastic} need (sometimes tight) bounds on the
accuracy of the function value proxy to achieve convergence both in
theory and in practice. They typically rely on second-order
information or higher, and achieve a much more favorable
$\mathcal{O}\left(\epsilon^{-3/2}\right)$ worst-case evaluation
complexity.

To achieve faster OFFO algorithms, the use of high-order derivatives
can therefore be considered. It is well-established that second-order methods
offer stronger theoretical guarantees than first-order methods, either
by demonstrating superior complexity for specific variants or by
being less sensitive to the problem’s conditioning in practice. This
approach has been successfully applied in proposing deterministic OFFO
variants of adaptive regularization \cite{OFFO-ARp} and trust region
methods \cite{Grapiglia2022}. These methods, while using significantly
less information, still achieve the complexity rate of
$\mathcal{O}(\epsilon^{-3/2})$ akin to their standard counterparts
\cite{Cartis2022-od}. Because they are not affected by errors in the
function value, they avoid the need for tight bounds on its accuracy,
making the algorithms more robust to noise. This robustness has been
confirmed in the numerical experiments proposed in
\cite{OFFO-ARp}. However, to the best of our knowledge, no theoretical
framework exists for high-order OFFO in the stochastic setting. 

We propose a theoretical framework for stochastic OFFO adaptive
regularization methods as introduced in \cite{OFFO-ARp}. Specifically,
we present an expected error on the approximative tensors up to the
$p$th order, with conditions dependent on the length of the previous  
$m$ steps. This approach allows for greater tolerance compared to work
that controls error with only the current or previous step
\cite{kohler17a, Agafonovetal2023}. By combining these tensor
conditions with classical probability and numerical analysis tools, we
can extend the results from the deterministic case in
\cite{OFFO-ARp}. Since our conditions depend only on past steps, the
implementation remains straightforward and can be adapted to machine
learning in terms of sampling sizes. The relaxed error bounds yield
promising results when applying the second-order algorithm. Our method
remains stochastic throughout, unlike other works where a
deterministic behavior is adopted at the end
\cite{kohler17a,RoostaKhorasani18,Bollapragadaetal18}. 

The paper is organized as follows: After restating the algorithm of
\cite{OFFO-ARp} and situating our condition on the probabilistic
derivatives within the literature in Section~2, we develop the
complexity rate analysis in Section~3. In Section~4, we outline some
potential applications of our algorithm. Initial numerical findings of
our algorithm for specific machine learning (ML) problems are
presented in Section~5. Finally, conclusions and perspectives are
drawn in Section~6. 

\numsection{A Stochastic OFFO adaptive regularization algorithm}~\label{algo-s}
\subsection{Problem Formulation}
\noindent

We consider the problem of finding approximate minimizers of the
unconstrained nonconvex optimization problem 
\beqn{problem}
\min_{x\in \Re^n} f(x),
\eeqn
where $f$ is a sufficiently smooth function from $\Re^n$ into $\Re$.
As motivated in the introduction, our aim is to design an algorithm in
which the objective function value is never computed and inexact
derivatives may be used. Our approach is based on regularization
methods. In such methods,  a model of the objective function is built
by ``regularizing'' a truncated inexact Taylor expansion of degree
$p$.  \\ 

\noindent We now detail the assumptions on the problems that we need
to establish our results.  

\begin{assumption}{assumption:AS1}
$f$ is $p$ times continuously differentiable in $\Re^n$.
\end{assumption}

\begin{assumption}{assumption:AS2}
There exists a constant $f_{\rm low}$ such that
$f(x) \geq f_{\rm low}$ for all $x \in \Re^n$.
\end{assumption}

\begin{assumption}{assumption:AS3}
The $p$th derivative of $f$ is globally Lipschitz continuous, that is,
there exists a non-negative constant $L_p$ such that 
\beqn{LipHessian}
\|\nabla_x^p f(x) - \nabla_x^p f(y)\| \leq L_p \|x-y\| \, \text{ for all } x, y \in \Re^n, \ms \text{ with } L_p \geq 3,
\eeqn
where $\| . \|$ denotes the Euclidean norm for vectors in $\Re^n$ and
$\| . \|$ the associated subordinate norm for $p$th order tensors.  
In the rest of the paper, all probabilistic approximations of exact
quantities will be denoted by an overline. 
\end{assumption}

\begin{assumption}{assumption:AS4}
If $p > 1$, there exists a constant $\kappa_{\text{high}} \geq 0$ such that
\beqn{negcurv}
\min_{\|d\|\leq 1} \overline{\nabla_x^i f}(x)[d]^i \geq -\kappa_{\text{high}} \text{ for all } x \in \Re^n \text{ et } i \in \{2, \ldots, p\},
\eeqn
where $\overline{\nabla_x^i f}(x) $ is the $i$th approximate
stochastic derivative tensor of $f$ computed at $x$ and where $T[d]^i$
denotes the $i$-dimensional tensor $T$ applied on $i$ copies of the
vector $d$. (For notational convenience, we set $\kappa_{\text{high}}
= 0$ if $p=1$).  
\end{assumption}

The previous assumptions \ref{assumption:AS1}--\ref{assumption:AS3}
are standard when studying the complexity of deterministic $p$th order
methods. Note that assumption \ref{assumption:AS4} is weaker than
imposing uniform boundedness on the sampled derivatives and is
standard in the study of Objective Free Function algorithms
\cite{OFFO-ARp}.  Moreover, it automatically holds in the exact case
for any function satisfying Assumption~\ref{assumption:AS1} on a
bounded domain. When $p=2$, the class of functions satisfying it is
sometimes (misleadingly) called "weakly convexity" and was shown to
cover different cases of interest, see both
\cite{DavisDrus19,Drusvyatskiy18} for more discussion on weak
convexity.

We will return to the probabilistic bounds that must
be satisfied by the tensor derivatives error after stating the
algorithm.

\subsection{The OFFO algorithm with stochastic derivatives}

Adaptive regularization methods are iterative schemes which compute a step from
an iterate $x_k$ to the next by approximately minimizing a $p$th
degree regularized model $m_k(s)$ of $f(x_k+s)$ of the form
\[ 
T_{f,p}(x_k,s) + \frac{\sigma_k}{(p+1)!} \|s\|^{p+1},
\]
where $T_{f,p}(x,s)$ is the $p$th order Taylor expansion of
functional $f$ at $x$ truncated at order $p$, that is,
\beqn{taylor}
T_{f,p}(x,s) \eqdef f(x) + \sum_{i=1}^p \frac{1}{i!} {\nabla_x^i f(x)}[d]^i. 
\eeqn
In particular, \ref{assumption:AS3} implies \cite[Corollary~A.8.4]{CartGoulToin22} that
\beqn{grad-lip}
\|\nabla_x^1f(x+s)-\nabla_x^1T_{f,p}(x,s)\|\leq\frac{L}{p!}\|s\|^p.
\eeqn
In the case where approximative  derivatives are used, one then uses an approximate $p$th order Taylor model 
\beqn{approxtaylor}
\overline{T_{f,p}}(x,s) \eqdef f(x) + \sum_{i=1}^p \frac{1}{i!} \overline{\nabla_x^i f}(x)[d]^i,
\eeqn
and the model $m_k$ is then,
\beqn{model}
m_k(s) \eqdef \overline{T_{f,p}}(x_k,s) + \frac{\sigma_k}{(p+1)!} \|s\|^{p+1}.
\eeqn
In \req{model}, the approximate $p$th order Taylor series is ``regularized'' by adding the term
$\frac{\sigma_k}{(p+1)!} \|s\|^{p+1}$, where $\sigma_k$ is known as the ``regularization
parameter''. This term guarantees that $m_k(s)$ is bounded below and thus
makes the procedure of finding a step $s_k$ by (approximately) minimizing
$m_k(s)$ well-defined. Our proposed algorithm follows the outline  of existing
\al{AR$p$} regularization methods
\cite{CartisGT11a,Birgin2016,Cartis2022-od} and the recent work of
\cite{OFFO-ARp} on an optimal $p$th order objective free function
method.  

We stress that unlike inexact adaptive second-order methods analyzed
in \cite{Yao2021, kohler17a, Bellavia2019, BellaviaGianBen2020}, we
don't evaluate the true function value nor a proxy.  In what follows,
all random quantities are denoted by capital letters, while the use of
small letters is reserved for their realization.

\algo{InexOFFARp}{ Stochastic  OFFO adaptive regularization of degree $p$ (\tal{StOFFAR$p$})}{
	\begin{description}
	   \item[Step 0: Initialization: ] An initial point $x_0\in \Re^n$, a regularization
		parameter $\sigma_0>0$ and a requested final gradient accuracy
		$\epsilon_1 \in (0,1]$ are given, as well as the parameter $\theta_1 > 1$.
		Set $k=0$.
	   \item[ Step 1: Compute current derivatives ] Evaluate
                $\overline{g_k} \eqdef \overline{\nabla_x^1 f}(x_k)$ and $\{\overline{\nabla_x^i f}(x_k) \}_{i=2}^p$. 
	      \item[Step 2: Step calculation: ]
                Compute a step $s_k$  which sufficiently reduces the
                model $m_k$ defined in \req{model} in the sense that 
		\beqn{descentmodel}
		 m_k(s_k) - m_k(0) \leq  0
		\eeqn
		and
		\beqn{gradstep}
                \|\overline{\nabla_s^1 T_{f,p}}(x_k,s_k)\| \leq  \theta_1 \frac{\sigma_k}{p!}\|s_k\|^p.
		\eeqn
	   \item[Step 3: Updates. ]
		Set
                \beqn{accept}
                x_{k+1} = x_k + s_k
                \eeqn
                and
		\beqn{vkupdate}
		\sigma_{k+1} =  \sigma_k + \sigma_k \| s_k\|^{p+1}.
		\eeqn
		Increment $k$ by one and go to Step~1.
	\end{description}
}
We now define the probabilistic notation which will be used throughout the paper. 
We emphasize that the approximate derivatives (as evaluated in Step~1)
are noisy random evaluations of the exact quantities. The
\al{StOFFAR$p$} algorithm therefore generates a stochastic process
\[
\{X_k, \, \overline{\nabla_x^i f}(X_k), \, \Sigma_k, \, S_k \}
\]
on some probability space ($\Omega$, $\calF$, $\mathbb{P}$). The
associated expectation operator will be denoted by $\E{.}$ and
$\Ek{.}{k}$ will stand for the conditional expectation  knowing 
$\left\{ \overline{\nabla_x^i f}(X_j) | j \in \iiz{k-1} \right\}$.
Note that $\Sigma_0  = \sigma_0$ is deterministic and we allow the
initialization $x_0$ to be a random variable. We also denote by $G_k
\eqdef \nabla_x^1 f(X_k) $ and
$\overline{G}_k \eqdef \overline{\nabla_x^1 f}(X_k)$.

\noindent 
One could have chosen $\sigma_{k}$ to be of the form
\beqn{practicalsigk}
\sigma_{k} \in [\vartheta \nu_k , \max (\nu_k, \eta_k) ] 
\eeqn
where $\eta_k$ is bounded non-negative sequence, $\nu_k$ is given by
\eqref{vkupdate} and $\vartheta$ a hyperparameter in $(0,1]$. This may
  be useful when devising a numerical implementation of the algorithm
  to closer adapt to local variations of the local Lipschitz constant,
  as done in \cite[Section 5]{OFFO-ARp}. We preferred to keep
  \eqref{vkupdate} for the sake of simplicity in our subsequent
  analysis. 

\noindent
Compared to the vanilla adaptive regularization methods where the
inequality in \eqref{descentmodel} must  be strict (see for example
\cite{Birgin2016}), we only require a simple decrease, since zero
derivatives can occur in the stochastic case.

The test \req{gradstep} follows from \cite{GratToin21} and extends the more
usual condition where the step $s_k$ is chosen to ensure that
\[
\|\nabla_s^1 m_k(s_k)\| \leq  \theta_1 \|s_k\|^p.
\]
It is indeed easy to verify that \req{gradstep} holds at a local
minimizer of $m_k$ with $\theta_1\geq 1$ (see  \cite{GratToin21} for
details).
Thus imposing \req{gradstep} and \req{descentmodel} simply
  amounts to minimze the model \req{model} inexactly. Dedicated
subroutines have been developed for the case $p=2$, see \cite[Chapter~8--10]{Cartis2022-od} and the references therein. For $p\geq3$, one can turn to a first-order algorithm (backtracking gradient descent \cite[Chapter~2]{Cartis2022-od}) or a second-order one such as standard trust-region \cite[Chapter~3]{Cartis2022-od}.

We propose the following conditions on the expectation of the
errors on the derivatives tensors.

\noindent 
\begin{assumption}{assumption:AS5}
There exists $\kappa_D \geq 0$ such that at each iteration $k \geq 0$, we have that 
\beqn{tensderiverror}
\Ek{\| \nabla_x^i f(X_k) - \overline{\nabla_x^i f}(X_k) \|^\sfrac{p+1}{p+1-i}}{k} \leq \kappa_D \xi_k, \tim{ for all } i \in \ii{p},
\eeqn 
\end{assumption}

\noindent
with $\xi_k = \sum_{i=1}^m\|S_{k-i}\|^{p+1}$ with the conventions that
\beqn{sigmajnegdef}
\| S_{-1} \| = \cdots \|S_{-m}\| \eqdef 1 \tim{ and }\sigma_{j} = \frac{\sigma_{0}}{2^{-j}}, \, \quad j \in \iibe{-m}{-1},
\eeqn
so that \eqref{vkupdate} is valid even when $k \in \iibe{-m}{-1}$.

We now discuss our proposed tensor conditions and compare them
with previously used requirements on stochastic derivatives, first focussing
on the case $m=1$ (we discuss the usage of \eqref{tensderiverror}
later).  Our subsequent discussion is divided into two parts: the first
considers different values of $p$, and the  second deals
with the practical case where $p=2$.
How to guarantee the conditions \eqref{tensderiverror} in
  practice will be discussed in Section~\ref{applicationsstooffo}.

First and foremost, note that the condition \eqref{tensderiverror} can
be related to the following  requirements on inexact tensors 
\beqn{tensderivcurrstep}
\| \nabla_x^i f(X_k) - \overline{\nabla_x^i f}(X_k) \| \leq \kappa_D \|S_{k} \|^{p-i+1}, \tim{ for all } i \in \ii{p},
\eeqn
 proposed both in \cite{Agafonovetal2023} and
 \cite[Chapter~13]{Cartis2022-od}.   However, one of the pitfalls of
 \eqref{tensderivcurrstep} is its implicit nature in that $S_{k}$ is
 not available when  $\overline{\nabla_x^i f}(X_k)$  is evaluated.
 Our condition \eqref{tensderiverror} only uses  information available
 from past iterations. Note that the exponent $\frac{p+1}{p+1-i}$  in \eqref{tensderiverror} is coherent with the standard condition \eqref{tensderivcurrstep} used to obtain the optimal complexity of tensor methods as proved in \cite[Chapter~13]{Cartis2022-od} or \cite{Agafonovetal2023}. This condition also implies that for a fixed $i$, a higher $p$ implies a tighter approximation of the $i$-th order tensor.
 
  To the best of the authors' knowledge, other
 stochastic adaptive regularization methods, such as that proposed by
 \cite{Bellavia2022}, additionally require more accurate function value approximations. 
  Specifically, the condition on the approximate function value
  $\overline{f}(x_k)$ in \cite{Bellavia2022} is that
\beqn{funcvaluecond}
|\overline{f}(x_k) - f(x_k)| \leq \eta \left(-\sum_{i=1}^p \frac{1}{i!} \overline{\nabla_x^i f}(x)[s_k]^i \right),
\eeqn
where $\eta$ is an algorithmic dependent constant and the  term in
parenthesis can be shown to be of order $\|s_k\|^{p+1}$, which is more restrictive than
\eqref{tensderivcurrstep}.  Moreover, the implicit bound
\req{funcvaluecond} must hold for all iterations, 
making it somewhat impractical for stochastic problems in machine learning
where subsampling is used. Note that the probabilistic
assumptions required for the approximate derivatives in 
\cite{Bellavia2022} do not treat each tensor derivative separately
but are slightly more general but also more abstract as they consider
their combined effect in the Taylor's expansion. For further details
on stochastic adaptive high-order methods, we refer the reader to
\cite{Bellavia2022} and the references therein.

We now turn to the case $ p = 2 $ and compare our framework with
previous stochastic cubic methods. The use of the past step to control
the error on the inexact gradient and the Hessian was first proposed
for numerical experiments in \cite{kohler17a} with good
empirical success, although the theory requires the use of the
current step as in \eqref{tensderivcurrstep}. This approach was
later investigated theoretically in \cite{WANG2019146} in an inexact
cubic regularization algorithm, where \eqref{tensderivcurrstep} is
assumed to hold with $\|S_{k-1}\|$ instead of $\|S_k\|$. However, the
authors unrealistically assume knowledge of the Lipschitz
constant. More recently, this idea has been combined with variance
reduction in \cite{Zhang22} to devise efficient cubic regularization
algorithms, but knowledge of the problem's geometry is still required. 

One drawback of using only the last step size to control the error is
that it may make the method exact after a few iterations, as
illustrated in \cite{kohler17a}. In contrast, we hope that using the
last $ m $ steps may provide better control. Intuitively, we are able
to use the last $ m $ steps to control the errors as our $\sigma_k$
update rule \eqref{vkupdate} accumulates the past steps size lengths. 

Other notable stochastic adaptive cubic regularization methods have
been developed in the literature; see for example
\cite{Tripuraneni18,ZhouXuGu19,BellaviaGurioli21,WangZhouYingbinLan19,ChaytiDoiJaggi23,ScheXie23,Zhang22}
to name a few. Let us now briefly review these references and
highlight the novelty of our approach. First, note that
\cite{Tripuraneni18,ZhouXuGu19,WangZhouYingbinLan19,ChaytiDoiJaggi23,Zhang22}
do not provide an adaptation mechanism for the regularization
parameter and typically assume knowledge of the Lipschitz Hessian
constant. We should also mention that the conditions proposed in 
\cite{ChaytiDoiJaggi23} are very similar to ours, as they also propose
bounds on 
\[
\Ek{\| \nabla_x^1 f(X_k) - \overline{\nabla_x^1 f}(X_k) \|^\frac{3}{2}}{k} \quad \text{and} \quad \Ek{\| \nabla_x^2 f(X_k) - \overline{\nabla_x^2 f}(X_k) \|^{3}}{k}.
\]
However, it should be noted that the analysis is
limited to the second-order case and again assumes the knowledge of
the Hessian Lipschitz constant. 

Another line of work \cite{BellaviaGurioli21,toint2020stochastic},
although adaptive, still requires an accurate approximation of the
function value to successfully adjust the regularization parameter
$\sigma_k$. Note that \cite{BellaviaGurioli21} also proposes inexact
conditions that are dynamic (as is the case here in
\eqref{tensderiverror}) and controlled by the inexact gradient norm.
Note also that this latter work imposes exact evaluation of the
objective-function value and is restricted to the second-order
case.  

Finally, a "fully" stochastic cubic method has been proposed in
\cite{ScheXie23}, where the gradient and the Hessian satisfy a
condition related to \eqref{tensderivcurrstep} with some
probability. However, they impose additional conditions on the
stochastic oracle of the function value. To our knowledge, no
practical case for machine learning has been proposed in
\cite{ScheXie23}. In contrast, our paper later proposes practical
variants of \eqref{tensderivcurrstep} for machine learning problems. 

\numsection{Evaluation complexity for the inexact \tal{StOFFAR$p$}
  algorithm}\label{complexity-s} 

\noindent
We start our analysis of evaluation complexity by defining
the following constants for notational convenience: 
\beqn{chi12def}
\chi_p^1 \eqdef \sum_{i=1}^p \frac{i}{i!(p+1)}, \quad \chi_p^2 \eqdef \sum_{i=1}^p \frac{p+i-1}{i!(p+1)}, \quad \kappa_p = \frac{(2p)^\sfrac{p+1}{p}}{(2p)} = (2p)^\sfrac{1}{p}.
\eeqn

\beqn{chi34def}
\chi_p^3 = \sum_{i=2}^p \left( \frac{(i-1)}{(i-1)! p} \right)^\sfrac{p+1}{p}, \quad \chi_p^4 = 1+ \sum_{i=2}^p \left( \frac{(p-i+1)}{(i-1)! p} \right)^\sfrac{p+1}{p}.
\eeqn

We may now state local decrease bounds resulting from classical Taylor
inequalities and the step computation mechanism. They combine the
standard bounds coming from AS.3 while also taking into account the
assumption on inexact derivatives AS.5 that holds in expectation.

\llem{lipschitz}{
  Suppose that AS.1, AS.3 and AS.5 hold and let $\alpha > 0$.  Then
  \beqn{Lip-f}
  \Ek{\frac{\Sigma_{k}}{(p+1)!} \|S_k\|^{p+1}}{k}
  \leq \Ek{f(X_{k}) - f(X_{k+1})  }{k} + \kap{a} \Ek{\|S_k\|^{p+1}}{k} + \kappa_D \chi_p^2 \xi_k
  \eeqn
  and
  \beqn{Lip-g}
   \Ek{\frac{\|G_{k+1}- \overline{\nabla_s^1 T_{f,p}}(X_k,S_k) \|^\sfrac{p+1}{p}}{\Sigma_{k+1}^\alpha}  }{k}
   \leq \kap{b} \Ek{\frac{\|S_k\|^{p+1}}{\Sigma_{k+1}^\alpha} }{k}
   + \kappa_p \kappa_D \chi_p^4 \frac{\xi_k}{\Sigma_{k}^\alpha},
  \eeqn
  where
  \beqn{ka-kb-def}
  \kap{a} \eqdef \frac{L_p}{(p+1)!} + \chi_p^1
  \tim{ and }
  \kap{b} \eqdef \kappa_p \left(\frac{L_p}{p!}\right)^\sfrac{p+1}{p} +  \chi_p^3.
  \end{equation}
 }

\proof{
By combining \eqref{taylor}, \eqref{approxtaylor}, \req{grad-lip} and basic tensor inequalities, we obtain that 
	 \begin{align*}
	 f(X_{k+1})- \overline{T_{f,p}}(X_k,S_k) &= 
	  f(X_{k+1}) - {T_{f,p}}(X_k,S_k) + {T_{f,p}}(X_k,S_k) - \overline{T_{f,p}}(X_k,S_k)  \\
	 &\leq  \frac{L_p}{(p+1)!} \|S_k\|^{p+1} + \sum_{i=1}^p \frac{1}{i!} \left| \left( \nabla_x^i f(S_k) - \overder{i}{k} \right) [S_k]^i \right| \\
	 &\leq \frac{L_p}{(p+1)!} \|S_k\|^{p+1} + \sum_{i=1}^p \frac{1}{i!} \|\nabla_x^i f(X_k) - \overder{i}{k} \| \| S_k\|^i.
	 \end{align*}
	 Using now Young's inequality with $p_i = \frac{p+1}{p+1-i}$ and
     $q_i = \frac{p+1}{i}$, we derive that
	 \begin{align*}
	 f(X_{k+1})- \overline{T_{f,p}}(X_k,S_k) \leq\ & \frac{L_p}{(p+1)!} \|S_k\|^{p+1} \\
	 &+  \sum_{i=1}^p \frac{1}{i!} \left( \frac{ (p+1-i)\|\nabla_x^i f(X_k) - \overder{i}{k} \|^\sfrac{p+1}{p+1-i}}{p+1} + \frac{i\|S_k \|^{p+1}}{p+1}  \right). 
	 \end{align*}
	 
	 Using the definition of $\chi_p^1$ in \eqref{chi12def}, that of
     $\kap{a}$ in \req{ka-kb-def} and the fact that
     \eqref{descentmodel} gives that
     $\overline{T_{f,p}}(X_k,S_k) \leq f(X_k) -
     \frac{\Sigma_{k}}{(p+1)!} \|S_k\|^{p+1}$, we obtain that
	 \begin{align*}
	  &\left(f(X_{k+1})- f(X_k) + \frac{\Sigma_{k}}{(p+1)!} \|S_k\|^{p+1} \right) \leq f(X_k) - \overline{T_{f,p}}(X_k,S_k)   \\ &\leq \kap{a} \|S_k\|^{p+1} +  \sum_{i=1}^p \frac{1}{i!} \left( \frac{ (p+1-i)\|\nabla_x^i f(X_k) - \overder{i}{k} \|^\sfrac{p+1}{p+1-i}}{p+1} \right).
	 \end{align*}
	 Taking now $\Ek{.}{k}$, using \eqref{tensderiverror} and rearranging
	 \[
	   \Ek{\frac{\Sigma_{k}}{(p+1)!} \|S_k\|^{p+1}}{k}
       \leq \Ek{f(X_{k+1}) - f(X_k)  }{k} + \kap{a}\Ek{\|S_k\|^{p+1}}{k}
	 + \kappa_D \sum_{i=1}^p \frac{p+1-i}{(p+1) i!}\, \xi_k. 
	 \]
	 Using now $\chi_p^2$ definition in \eqref{chi12def}, we obtain \eqref{Lip-f}.
	 
	 \noindent
	 We turn now to the proof of \eqref{Lip-g}. Using the triangle
     inequality, \req{grad-lip}, \eqref{taylor} and
     \eqref{approxtaylor} yields that
	 \begin{align*}
	 \|G_{k+1}- \overline{\nabla_s^1 T_{f,p}}(X_k,S_k) \| &\leq \|G_{k+1}- {\nabla_s^1 T_{f,p}}(X_k,S_k) \| + \|\nabla_s^1 T_{f,p}(X_k,S_k)- \overline{\nabla_s^1 T_{f,p}}(X_k,S_k) \| \\
 &\leq  \frac{L_p}{p!} \|S_k\|^{p} + \sum_{i=1}^p \frac{1}{(i-1)!} \left| \left( \nabla_x^i f(X_k) - \overline{\nabla_x^i f}(X_k) \right) [S_k]^{i-1} \right| \\
 &\leq \frac{L_p}{p!} \|S_k\|^{p} + \sum_{i=1}^p \frac{1}{(i-1)!} \|\nabla_x^i f(X_k) - \overline{\nabla_x^i f}(X_k) \| \| S_k\|^{i-1} \\
 &\leq \frac{L_p}{p!} \|S_k\|^{p} + \|\nabla_x^1 f(X_k) - \overline{\nabla_x^1 f}(X_k) \|  \\ 
 &\hspace*{1em}+  \sum_{i=2}^p \frac{1}{(i-1)!} \|\nabla_x^i f(X_k) - \overline{\nabla_x^i f}(X_k) \| \| S_k\|^{i-1}.
     \end{align*}
Again using Young's inequality with $p_i = \frac{p}{p+1-i}$ and $q_i =
\frac{p}{i-1}$ for $i \in \iibe{2}{p}$, we derive that
\begin{align*}
\|G_{k+1}- \overline{\nabla_s^1 T_{f,p}}(X_k,S_k) \| \leq\ & \frac{L_p}{p!} \|S_k\|^{p} + \|\nabla_x^1 f(X_k) - \overline{\nabla_x^1 f}(X_k) \| \\  &+  \sum_{i=2}^p \frac{1}{(i-1)!} \left( \frac{p+1-i}{p} \|\nabla_x^i f(X_k) - \overline{\nabla_x^i f}(X_k) \|^\sfrac{p}{p+1-i} + \frac{i-1}{p} \| S_k\|^{p} \right).
\end{align*}
Taking the last inequality to the $\frac{p+1}{p}$ power, using the
fact that $x^\frac{p+1}{p}$ is a convex function, the definition of
$\kappa_p$ in \eqref{chi12def}, the fact that the left-hand side has $2p$ terms
and dividing both sides of the inequality by $\Sigma_{k+1}^{\alpha}$
gives that
\begin{align*}
\frac{\|G_{k+1}- \overline{\nabla_s^1 T_{f,p}}(X_k,S_k)
  \|^\sfrac{p+1}{p}}{\Sigma_{k+1}^{\alpha}} \leq\ & \kappa_p \left(\frac{L_p}{p!} \right)^\sfrac{p+1}{p} \frac{\|S_k\|^{p+1}}{\Sigma_{k+1}^{\alpha} } + \kappa_p \frac{\|\nabla_x^1 f(X_k) - \overline{\nabla_x^1 f}(X_k) \|^\sfrac{p+1}{p}}{\Sigma_{k+1}^{\alpha}} \\  &+ \kappa_p  \sum_{i=2}^p \frac{1}{(i-1)!^\sfrac{p+1}{p}} \Biggl( \left[\frac{p+1-i}{p} \right]^\sfrac{p+1}{p} \frac{\|\nabla_x^i f(X_k) - \overline{\nabla_x^i f}(X_k) \|^\sfrac{p+1}{p+1-i}}{\Sigma_{k+1}^{\alpha}} \\ &+ \left[\frac{i-1}{p} \right]^\sfrac{p+1}{p} \frac{\| S_k\|^{p+1}}{\Sigma_{k+1}^{\alpha}} \Biggr).
\end{align*}

Taking the conditional expectation over the past iterations, using
\eqref{tensderiverror} and the fact that
$\frac{1}{\Sigma_{k+1}^\alpha} \leq \frac{1}{\Sigma_{k}^\alpha}$ for
the terms $\|\nabla_x^i f(X_k) - \overline{\nabla_x^i f}(X_k) \|$, we
derive that 
\begin{align*}
\Ek{ \frac{\|G_{k+1}- \overline{\nabla_s^1 T_{f,p}}(X_k,S_k)
    \|^\sfrac{p+1}{p}}{\Sigma_{k+1}^\alpha} }{k} \leq\ & \kappa_p \left(\frac{L_p}{p!}\right)^\sfrac{p+1}{p} \Ek{ \frac{\|S_k \|^{p+1}}{\Sigma_{k+1}^\alpha} }{k} + \kappa_p \kappa_D \frac{\xi_k}{\Sigma_{k}^\alpha} \\ &+ \kappa_p \sum_{i=2}^p \frac{1}{(i-1)!^\sfrac{p+1}{p}} \Biggl( \left[\frac{p+1-i}{p} \right]^\sfrac{p+1}{p} \kappa_D  \frac{\xi_k}{\Sigma_{k}^\alpha}\\ &+ \left[\frac{i-1}{p} \right]^\sfrac{p+1}{p} \Ek{ \frac{\| S_k\|^{p+1}}{\Sigma_{k+1}^\alpha} }{k} \Biggr). 
\end{align*}
Using now the definitions of $\chi_p^3$ and $\chi_p^4$ in
\eqref{chi34def} and that of $\kap{b}$ in \req{ka-kb-def}, we obtain (\ref{Lip-g}).
} 

\noindent

The next lemma provides two useful upper bounds on the gradient norm
at iteration $ k+1 $ divided by the regularization parameter. This
also clarifies why Lemma~\ref{lipschitz} was stated with a generic
$\alpha$ parameter: we will need this result for two different values
of $\alpha$ in the following proof, which is inspired by
\cite[Lemma~2.3]{Birgin2016}, but it also takes into account the
derivative tensor errors that hold in expectation (AS.5) and the
update rule of the regularization parameter $\sigma_{k}$ in
\eqref{vkupdate}. 

\llem{useful}{
	Suppose that AS.1, AS.3 and AS.5 hold and let $k \geq 0$. Then,
	 \beqn{crucial}
	 \Ek{\frac{\|{G_{k+1}}\|^\sfrac{p+1}{p}}{\Sigma_{k+1}^\sfrac{p+1}{p}}
     }{k} \leq  
     \kap{c}\Ek{\frac{\Sigma_{k+1}-\Sigma_{k}}{\Sigma_{k+1}} } {k}
	 + \kap{d} \sum_{j=k-m}^{k-1}\frac{\Sigma_{j+1}-\Sigma_{j}}{\Sigma_{j+1}}
	 \eeqn
     and
\beqn{crucial2}
\Ek{\frac{\|{G_{k+1}}\|^\sfrac{p+1}{p}}{\Sigma_{k+1}^\sfrac{1}{p}}}{k}
\leq  \kap{c}  \Ek{\Sigma_{k+1}-\Sigma_{k}} {k} + \kap{d} \sum_{j=k-m}^{k-1}[\Sigma_{j+1}-\Sigma_{j}],
\eeqn
where 
\beqn{kc-kd-def}
\kap{c} \eqdef\frac{\twoppp}{\sigma_{0}^\sfrac{p+1}{p}}\left( \kap{b} + \frac{\theta_1^\sfrac{p+1}{p}}{p!^\sfrac{p+1}{p}}\sigma_{0}^\sfrac{p+1}{p}\right)
\tim{ and }
\kap{d} \eqdef \frac{2^\sfrac{mp+1}{p} \kappa_p \kappa_D \chi_p^4}{\sigma_{0}^\sfrac{p+1}{p}},
\eeqn
with $\kappa_p$, $\chi_p^3$ and $\chi_p^4$  defined in
\eqref{chi12def} and \eqref{chi34def} and $\kap{b}$ given by \req{ka-kb-def}.     
}
\proof{
First consider $\alpha \in \{\frac{1}{p}, \, \frac{p+1}{p} \}$.
	From the triangular inequality and the fact that
    $(x+y)^\sfrac{p+1}{p} \leq \twoppp \left( x^\sfrac{p+1}{p} +
    y^\sfrac{p+1}{p} \right)$ for $x,y \geq 0$ and \eqref{gradstep},
    we obtain  
	\begin{align*}
		\frac{\|{G_{k+1}}\|^\sfrac{p+1}{p}}{\Sigma_{k+1}^\alpha} &\leq  \frac{\left(\|G_{k+1}- \overline{\nabla_s^1 T_{f,p}}(X_k,S_k) \| + \|\overline{\nabla_s^1 T_{f,p}}(X_k,S_k) \|\right)^\sfrac{p+1}{p}}{\Sigma_{k+1}^\alpha} \\
		&\leq   2^\sfrac{1}{p} \frac{\|G_{k+1}- \overline{\nabla_s^1 T_{f,p}}(X_k,S_k) \|^\sfrac{p+1}{p}}{\Sigma_{k+1}^\alpha} + 2^\sfrac{1}{p} \frac{\|\overline{\nabla_s^1 T_{f,p}}(X_k,S_k) \|^\sfrac{p+1}{p}}{\Sigma_{k+1}^\alpha} \\
		&\leq \twoppp \left( \frac{\|G_{k+1}- \overline{\nabla_s^1 T_{f,p}}(X_k,S_k) \|^\sfrac{p+1}{p}}{\Sigma_{k+1}^\alpha}  +  \frac{ (\theta_1\Sigma_{k})^\sfrac{p+1}{p}}{p!^\sfrac{p+1}{p}\Sigma_{k+1}^\alpha} \|S_k \|^{p+1}\right).
	\end{align*}
Taking $\Ek{.}{k}$ and using \eqref{Lip-g}, we derive that
	\begin{align}\label{twoalphas}
\Ek{\frac{\|{G_{k+1}}\|^\sfrac{p+1}{p}}{\Sigma_{k+1}^\alpha} }{k} 
\leq\ & \twoppp  \kap{b} \Ek{\frac{\|S_k\|^{p+1}}{\Sigma_{k+1}^\alpha}}{k}
+ \twoppp  \kappa_p \kappa_D \chi_p^4 \frac{\xi_k}{\Sigma_{k}^\alpha} 
+ \twoppp  \frac{\theta_1^\sfrac{p+1}{p}}{p!^\sfrac{p+1}{p}}
\Ek{\Sigma_{k}^\sfrac{p+1}{p}\frac{\|S_k\|^{p+1}}{\Sigma_{k+1}^\alpha}
}{k}.
\end{align}		
	
We first prove \eqref{crucial} and start with $\alpha = \frac{p+1}{p}$.
Using that $\xi_k = \sum_{j= 1}^{m} \|S_{k-j}\|^{p+1}$, that
$\|S_j\|^{p+1} = \frac{\Sigma_{j+1} - \Sigma_{j}}{\Sigma_{j}}$ for $j
\in \iibe{k-m}{k}$, \eqref{twoalphas}, and also that
$\Sigma_{k}$ is non-decreasing, also $\Sigma_{k} \geq \sigma_0$ ,
$\Sigma_{k - m} \geq \frac{\sigma_0}{2^m}$ for $k \geq 0$,  both facts
resulting from \eqref{vkupdate} and \eqref{sigmajnegdef}, we derive that
\begin{align*}
\Ek{\frac{\|{G_{k+1}}\|^\sfrac{p+1}{p}}{\Sigma_{k+1}^\sfrac{p+1}{p}} }{k} 
\leq\ & \twoppp  \kap{b} \Ek{\frac{\Sigma_{k+1}-\Sigma_{k}}{\Sigma_{k+1}} \frac{1}{\Sigma_{k} \Sigma_{k+1}^\sfrac{1}{p}} } {k}
+ \twoppp \kappa_p \kappa_D \chi_p^4 \sum_{j=k-m}^{k-1} 
\frac{\Sigma_{j+1}-\Sigma_{j}}{\Sigma_{j}} \frac{1}{\Sigma_{j+1}
  \Sigma_{k}^\sfrac{1}{p}}\\
& + \twoppp \frac{\theta_1^\sfrac{p+1}{p}}{p!^\sfrac{p+1}{p}} \Ek{\frac{\Sigma_{k+1}-\Sigma_{k}}{\Sigma_{k+1}} \frac{\Sigma_{k}^\sfrac{1}{p}}{\Sigma_{k+1}^\sfrac{1}{p}} }{k} \\
\leq\ & \frac{\twoppp}{\sigma_{0}^\sfrac{p+1}{p}}  \kap{b} \Ek{\frac{\Sigma_{k+1}-\Sigma_{k}}{\Sigma_{k+1}} } {k}
+ \kap{d}  \sum_{j=k-m}^{k-1} \frac{\Sigma_{j+1}-\Sigma_{j}}{\Sigma_{j+1}} 
+ \twoppp \frac{\theta_1^\sfrac{p+1}{p}}{p!^\sfrac{p+1}{p}} \Ek{\frac{\Sigma_{k+1}-\Sigma_{k}}{\Sigma_{k+1}} }{k},
\end{align*}
where $\kap{d}$ is defined in \req{kc-kd-def}.
Rearranging the last inequality and using the definition of $\kap{c}$
in \req{kc-kd-def} yields inequality \eqref{crucial}.

Consider now the case where $\alpha = \frac{1}{p}$.
Again, using the same arguments used to prove \eqref{crucial}, we deduce that
\begin{align*}
\Ek{\frac{\|{G_{k+1}}\|^\sfrac{p+1}{p}}{\Sigma_{k+1}^\sfrac{1}{p}} }{k} 
\leq\ & \twoppp  \kap{b} \Ek{\frac{\Sigma_{k+1}-\Sigma_{k}}{\Sigma_{k}} \frac{1}{\Sigma_{k+1}^\sfrac{1}{p}} } {k}
+ \twoppp \kappa_p \kappa_D \chi_p^4 \sum_{j=k-m}^{k-1} 
\frac{\Sigma_{j+1}-\Sigma_{j}}{\Sigma_{j}} \frac{1}{\Sigma_{k}^\sfrac{1}{p}} \\  &+ \twoppp  \frac{\theta_1^\sfrac{p+1}{p}}{p!^\sfrac{p+1}{p}} \Ek{ \frac{\Sigma_{k}^\sfrac{1}{p} \left( \Sigma_{k+1}-\Sigma_{k}\right)}{\Sigma_{k+1}^\sfrac{1}{p}} }{k} \\
\leq\ & \frac{\twoppp}{\sigma_{0}^\sfrac{p+1}{p}}  \kap{b} \Ek{\Sigma_{k+1}-\Sigma_{k} } {k}
+ \kap{d} \sum_{j=k-m}^{k-1} [\Sigma_{j+1}-\Sigma_{j}]
+ \twoppp \frac{\theta_1^\sfrac{p+1}{p}}{p!^\sfrac{p+1}{p}} \Ek{\Sigma_{k+1}-\Sigma_{k} }{k}.
\end{align*}
Rearranging the last inequality gives the second result of the lemma.
} 

The following lemma restates a result similar to that developed
when analyzing the exact version of Algorithm~\ref{InexOFFARp} in
\cite{OFFO-ARp}, but we extend it by providing a bound on
$\|S_k\|^{p+1}$, under the assumption that $\|S_k\|^{p+1}$ is bounded
by a constant depending on AS.4. 

\llem{stepgkbound}{
Suppose that AS.1 and AS.4 hold. At each iteration $k$, we have that
\beqn{skbound}
\| S_k\| \leq 2 \max\left( \eta,  \left(\frac{(p+1)! \| \overline{G_k} \|}{\Sigma_k}\right)^\sfrac{1}{p} \right),
\eeqn
where 
\beqn{kaphigh}
\eta = \max_{i \in \iibe{2}{p}}  \left[\frac{\khigh(p+1)!}{i!\, \sigma_0}\right]^\sfrac{1}{p-i+1}.
\eeqn
Moreover, 
\beqn{SKboundizycase}
\|S_k\|^{p+1} \indica{\|S_k\| \leq 2\eta} \leq (1+2^{p+1} \eta^{p+1} ) \frac{\Sigma_{k+1} - \Sigma_{k}}{\Sigma_{k+1}} \indica{\|S_k\| \leq 2\eta}. 
\eeqn
}
\proof{See Appendix~\ref{skborne} for the proof of inequality \eqref{skbound}.
We now turn to establishing \eqref{SKboundizycase}. From
\eqref{vkupdate}, and the fact that $\|S_k\|^{p+1} \indica{\|S_k\|
  \leq 2\eta} \leq (2\eta)^{p+1}$, we have that
\[
\Sigma_{k+1} \indica{\|S_k\| \leq 2\eta} = \Sigma_{k} \indica{\|S_k\| \leq 2\eta} + \|S_k\|^{p+1} \Sigma_{k } \indica{\|S_k\| \leq 2\eta} \leq \Sigma_{k} \indica{\|S_k\| \leq 2\eta} \left( 1 + (2\eta)^{p+1} \right),
\]
which yields that
\[
\frac{\Sigma_{k+1} \indica{\|S_k\| \leq 2\eta}}{1 + (2\eta)^{p+1}} \leq \Sigma_{k} \indica{\|S_k\| \leq 2\eta}.
\]

Multiplying both sides of the previous inequality by $\|S_k\|^{p+1}$,
adding $\Sigma_{k} \indica{\|S_k\| \leq 2\eta}$, and  using identity
\eqref{vkupdate}, we derive that 
\begin{align*}
\Sigma_{k} \indica{\|S_k\| \leq 2\eta} + \frac{\Sigma_{k+1}}{1+(2\eta)^{p+1}} \|S_k\|^{p+1} \indica{\|S_k\| \leq 2\eta}   \leq \indica{\|S_k\| \leq 2\eta} \left( \Sigma_{k} + \Sigma_{k}  \|S_k\|^{p+1} \right) = 
\Sigma_{k+1} \indica{\|S_k\| \leq 2\eta}.
\end{align*}
Now rearranging the last inequality  yields \eqref{SKboundizycase}.
} 

Why we have showcased the term
$\frac{\Sigma_{k+1}-\Sigma_{k}}{\Sigma_{k+1}}$ in the upper bound of
\eqref{SKboundizycase}, it will become clear later in the paper. We
now turn to proving a bound on $\E{\left( \frac{\|\overline{G_k}
    \|}{\Sigma_k}\right)^\sfrac{p+1}{p}}$, as it will allow us to
derive a bound on $\E{\|S_k\|^{p+1}}$. 

\llem{Gkoverinexcpec}{Suppose that AS.1, AS.3 and AS.5 hold and
  consider an iteration $k \geq 1$. Then, we have that 
	\begin{align}\label{Gkbound}
	  \E{\left(\frac{\|\overline{G_k}\|}{\Sigma_k}\right)^\sfrac{p+1}{p}}
      \leq\ & \frac{\kappa_D 2^\sfrac{mp+2}{p}}{\sigma_{0}^\sfrac{p+1}{p} }\hspace*{-1.5mm} \sum_{j=k-m}^{k-1}\E{\frac{\Sigma_{j+1} - \Sigma_{j}}{\Sigma_{j+1}}} +  2^\sfrac{1}{p}\kap{d} \sum_{j=k-m-1}^{k-2}\hspace*{-1mm}\E{\frac{\Sigma_{j+1}-\Sigma_{j}}{\Sigma_{j+1}}} \nonumber
	 \\
	&+ 2^\sfrac{1}{p} \kap{c} \E{\frac{\Sigma_{k} - \Sigma_{k - 1}}{\Sigma_{k}}},
	\end{align}
    where $\kap{c}$ is given by \req{kc-kd-def}.
We also have that, for $k = 0$,
\beqn{overG0bound}
\E{\left(\frac{\|\overline{G_0}\|}{\sigma_0}\right)^\sfrac{p+1}{p}} \leq \twoppp \left(\frac{\kappa_D}{\sigma_0^\sfrac{p+1}{p}} + \frac{\E{\|G_0\|^{\sfrac{p+1}{p}}}}{\sigma_{0}^\sfrac{p+1}{p}} \right).  
\eeqn 
}
\proof{
Consider an arbitrary positive $k$. From the triangle inequality and
the fact that $(x+y)^\sfrac{p+1}{p} \leq \twoppp \left(
x^\sfrac{p+1}{p} + y^\sfrac{p+1}{p} \right)$ for $x,y \geq 0$, we
derive that 
\beqn{decomposition}
\left(\frac{\|\overline{G_k}\|}{\Sigma_k}\right)^\sfrac{p+1}{p} 
\leq \left(\frac{\|\overline{G_k} - G_k\| + \|G_k\|}{\Sigma_k}\right)^\sfrac{p+1}{p} \leq \twoppp \left( \frac{\|\overline{G_k} - G_k\|^\sfrac{p+1}{p}}{\Sigma_k^\sfrac{p+1}{p}} + \frac{\|G_k\|^\sfrac{p+1}{p}}{\Sigma_k^\sfrac{p+1}{p}} \right).
\eeqn
Taking $k=0$, using the fact that $\Sigma_{k } = \sigma_0$ and
\eqref{tensderiverror} with $i=1$ yields \eqref{overG0bound}. 

Consider now $k \geq 1$. From the inequality \eqref{decomposition}, it
is sufficient to provide a bound on the two terms of the left-hand
side in order to establish the lemma's result. 	

Let us first provide a bound on
$\frac{\|\overline{G_k}-G_k\|^\sfrac{p+1}{p}}{\Sigma_k^\sfrac{p+1}{p}}$ in expectation. 
Using that $\Sigma_{k}$ is measurable with respect to the past, 
\eqref{tensderiverror} for $i=1$, \eqref{vkupdate} and that
$\Sigma_{j} \geq \frac{\sigma_{0}}{2^m}$ for $j \geq -m$ from
\eqref{sigmajnegdef} and $\Sigma_{k} \geq \sigma_{0}$ for $k \geq 0$
and that $\xi_k = \sum_{j= 1}^m \|S_{k-j}\|^{p+1}$, we derive  
\begin{align*}
	\Ek{\frac{\|\overline{G_k} - G_k\|^\sfrac{p+1}{p}}{\Sigma_k^\sfrac{p+1}{p}}}{k} &= \frac{1}{\Sigma_k^\sfrac{p+1}{p}} \Ek{\|\overline{G_k} - G_k\|^\sfrac{p+1}{p} }{k} \\
	&\leq \frac{\kappa_D}{\Sigma_k^\sfrac{p+1}{p}} \xi_k =\frac{\kappa_D}{\Sigma_k^\sfrac{p+1}{p}} \sum_{j=k-m}^{k-1}\frac{\Sigma_{j+1} - \Sigma_{j}}{\Sigma_{j}} 
	\leq \frac{\kappa_D 2^{\sfrac{mp+1}{p}}}{\sigma_{0}^\sfrac{p+1}{p}} \sum_{j=k-m}^{k-1}\frac{\Sigma_{j+1} - \Sigma_{j}}{\Sigma_{j+1}} \frac{}{}. 
\end{align*} 	
Taking the full expectation in the last inequality yields
\beqn{premierterme}
\E{\frac{\|\overline{G_k} - G_k\|^\sfrac{p+1}{p}}{\Sigma_k^\sfrac{p+1}{p}}} \leq \frac{\kappa_D 2^\sfrac{mp+1}{p}}{\sigma_{0}^\sfrac{p+1}{p}} \sum_{j=k-m}^{k-1} \E{ \frac{\Sigma_{j+1} - \Sigma_{j}}{\Sigma_{j+1}}}. 
\eeqn
Let us now focus on the second term
$\frac{\|G_k\|^\sfrac{p+1}{p}}{\Sigma_k^\sfrac{p+1}{p}}$. Using
\eqref{crucial}, and taking the full expectation yields that 
\beqn{secondterme}
\E{\frac{\|{G_{k}}\|^\sfrac{p+1}{p}}{\Sigma_{k}^\sfrac{p+1}{p}}  }
\leq  
 \kap{c}  \E{\frac{\Sigma_{k}-\Sigma_{k-1}}{\Sigma_{k}} }
+ \kap{d}  \sum_{j=k-m-1}^{k-2}\E{\frac{\Sigma_{j+1}-\Sigma_{j}}{\Sigma_{j+1}}}.
\eeqn
Note that as $k \geq 1$, the last term in the right-hand side of the
previous inequality is always well-defined. Injecting now
\eqref{premierterme} and \eqref{secondterme} in \eqref{decomposition}
gives the desired result. 
} 

The next lemma explains the specific dependency of the bounds
\eqref{Gkbound} and \eqref{SKboundizycase} with respect to 
$\Sigma_k$. 

\llem{ajsum}{Let $\{(a_j)_{j \in \iibe{m}{n}}\}$ be a positive
  nondecreasing sequence with $m < n$ and $(m,\,n) \in \mathbb{Z}^2$
  . Then, we have that 
\beqn{logsumbound}
\sum_{j=m+1}^n \frac{a_j - a_{j-1}}{a_j} \leq \log\left( a_n \right) - \log \left( a_{m}\right). 
\eeqn
}
\proof{
Let $j \geq m+1$ and suppose that $a_{j} > a_{j-1}$. By using the
concavity of the logarithm and since the sequence $a_i$ is
positive and further rearranging, we derive  
\[
\frac{a_j - a_{j-1}}{a_{j}} \leq \log(a_j) - \log(a_{j-1}).
\]
Note that the last inequality is still valid even when $a_j = a_{j-1}$. Thus, summing the last inequality for $j \in \iibe{m+1}{n}$ yields \eqref{logsumbound}.
} 

Combining the results of Lemmas~\ref{stepgkbound},~\ref{Gkoverinexcpec} and \ref{ajsum}, we are
now able to provide a bound on $\sum_{j= 0}^k \E{\|S_j\|^{p+1}}$. 

\llem{sigmakbound}{Suppose that AS.1, AS.3, AS.4 and AS.5 hold.
	 Then
\beqn{Esumstep}
\sum_{j= 0}^k \E{\|S_j\|^{p+1}} \leq \kap{e} + \kap{f} \log(\E{\Sigma_{k+1}}), 
\eeqn
where $\kap{e}$ and $\kap{f}$ are defined by
{ \small
\begin{align}\label{ke-def}
\kap{e} &\eqdef 2^{p+1} \pplusonefacpow \Biggl( \sumskzeroterm + 2^\sfrac{1}{p}\kap{d} \left(\log(2)\frac{m+1}{2}-\log(\sigma_{0})\right) \nonumber \\
&-2^\sfrac{1}{p} \kap{c}\log(\sigma_{0}) + \frac{\kappa_D 2^\sfrac{mp+2}{p}m}{\sigma_{0}^\sfrac{p+1}{p} } \left( \log(2) \frac{m-1}{2} -\log(\sigma_{0}) \right) \Biggr) 
- (1+2^{p+1} \eta^{p+1}) \log(\sigma_{0}),
\end{align}
where $\kap{c}$ is given by \req{kc-kd-def} and
\beqn{kf-def}
\kap{f} \eqdef 1+2^{p+1} \eta^{p+1}
+  \frac{2^{p+1} \pplusonefacpow}{\sigma_{0}^\sfrac{p+1}{p}}
\left(\kappa_D 2^\sfrac{mp+2}{p}m +  2^\sfrac{mp+2}{p} \chi_p^4 \kappa_D \kappa_pm +  2^\sfrac{2}{p} \kap{c}\right).
\eeqn
}
}
\proof{
From inequalities \eqref{skbound} and \eqref{SKboundizycase}, we have that

\begin{align*}
\|S_j\|^{p+1} 
&\leq \|S_j\|^{p+1} \indica{\|S_j\| \leq 2 \eta} + \|S_j\|^{p+1} \indica{\|S_j\|  \leq 2 \left( \frac{(p+1)! \| \overline{G}_j\|}{\Sigma_j}\right)^\sfrac{1}{p}} \nonumber \\ 
&\leq (1+2^{p+1} \eta^{p+1}) \frac{\Sigma_{j+1} - \Sigma_{j}}{\Sigma_{j+1}} +  2^{p+1} \left(  \frac{(p+1)! \| \overline{G}_j\|}{\Sigma_j}\right)^\sfrac{p+1}{p}.
\end{align*}
 Taking the full expectation in the last inequality, summing for $j
 \in \iiz{k}$, using \eqref{Gkbound} for $j\geq1$ and
 \eqref{overG0bound} when $j=0$, we derive that
 \begin{align*}
 \sum_{j= 0}^k \E{\|S_j\|^{p+1}} 
 &\leq (1+2^{p+1} \eta^{p+1}) \sum_{j= 0}^k \E{\frac{\Sigma_{j+1} -
     \Sigma_{j}}{\Sigma_{j+1}}} + 2^{p+1} \pplusonefacpow  \twoppp
 \left(\frac{\kappa_D}{\sigma_0^\sfrac{p+1}{p}}
 + \frac{\E{\|G_0\|^{\sfrac{p+1}{p}}}}{\sigma_{0}^\sfrac{p+1}{p}} \right) \\  
 &+ 2^{p+1} \pplusonefacpow \bigsum_{j=1}^k \Biggl(
 \frac{\kappa_D 2^\sfrac{mp+2}{p}}{\sigma_{0}^\sfrac{p+1}{p} }\left[  \sum_{\ell = j-m}^{j-1}\E{\frac{\Sigma_{\ell+1} - \Sigma_{\ell}}{\Sigma_{\ell+1}}} +    
\chi_p^4 \kappa_p
 \sum_{\ell = j-m-1}^{j-2}\E{\frac{\Sigma_{\ell+1} - \Sigma_{\ell}}{\Sigma_{\ell+1}}} \right] \nonumber\\
 &+ \frac{2^\sfrac{2}{p}}{\sigma_{0}^\sfrac{p+1}{p}} \kap{c} \E{\frac{\Sigma_{j} - \Sigma_{j - 1}}{\Sigma_{j}}}  \Biggr).
 \end{align*}
 We now provide a bound on the two sums involving
 $\E{\frac{\Sigma_{\ell+1} -\Sigma_{\ell}}{\Sigma_{\ell+1}}}$.
 By inverting the two sums, the linearity of the expectation,
 Lemma~\ref{ajsum} and the fact that $\Sigma_{k}$ is non-decreasing,
 we derive after some simplification,  that 
\begin{align}\label{firstboundellj}
\sum_{j= 1}^k \sum_{\ell = j-m}^{j-1} \E{\frac{\Sigma_{\ell+1} - \Sigma_{\ell}}{\Sigma_{\ell+1}}} &= \sum_{\ell = 1}^{m} \sum_{j= 1}^k  \E{\frac{\Sigma_{j-\ell+1} - \Sigma_{j-\ell}}{\Sigma_{j-\ell+1}}} \nonumber \\
&\leq  \sum_{\ell = 1}^{m} \E{\log(\Sigma_{k-\ell+1}) -
  \log(\Sigma_{1-\ell})} \leq m \E{\log({\Sigma_k})} - \sum_{\ell=1}^{m}
\log\left(\frac{\sigma_{0}}{2^{1-l}}\right) \nonumber \\ 
&\leq m (\E{\log(\Sigma_{k})} - \log(\sigma_{0})) + \log(2) \frac{m^2-m}{2}.
\end{align}
Similarly for $\sum_{j=1}^{k} \sum_{\ell=j-m-1}^{j-2}
\E{\frac{\Sigma_{\ell+1} - \Sigma_{\ell}}{\Sigma_{\ell+1}}}$ and using
the same arguments as  above yields that
\beqn{secboundellj}
\sum_{j= 1}^k \sum_{\ell = j-m-1}^{j-2} \E{\frac{\Sigma_{\ell+1} - \Sigma_{\ell}}{\Sigma_{\ell+1}}} = \sum_{\ell = 1}^{m} \sum_{j= 1}^k  \E{\frac{\Sigma_{j-\ell} - \Sigma_{j-\ell-1}}{\Sigma_{j-\ell}}} \leq m(\E{\log(\Sigma_{k-1})} - \log(\sigma_{0})) + \log(2) \frac{m^2+m}{2}.
\eeqn 
 
Now using \eqref{firstboundellj}, \eqref{secboundellj} and the
linearity of the expectation, we obtain that 
\begin{align*}
 \sum_{j= 0}^k \E{\|S_j\|^{p+1}} 
 &\leq (1+2^{p+1} \eta^{p+1})  \E{\log(\Sigma_{k+1}) - \log(\sigma_{0})} + 2^{p+1} \pplusonefacpow  \twoppp \left(\frac{\kappa_D}{\sigma_0^\sfrac{p+1}{p}} + \frac{\E{\|G_0\|^{\sfrac{p+1}{p}}}}{\sigma_{0}^\sfrac{p+1}{p}} \right) \\  
 &+ 2^{p+1} \pplusonefacpow  \Biggl( \frac{\kappa_D 2^\sfrac{mp+2}{p}}{\sigma_{0}^\sfrac{p+1}{p} } m\left( \E{\log(\Sigma_{k}) - \log(\sigma_{0})} + \frac{m-1}{2}\log(2)\right)\\ 
 &+    
 \frac{2^\sfrac{mp+2}{p} \chi_p^4 \kappa_D \kappa_p}{\sigma_{0}^\sfrac{p+1}{p}}m \left(  \E{\log(\Sigma_{k-1}) - \log(\sigma_{0})} + \log(2) \frac{m+1}{2} \right) \nonumber
 \\
 &+ \frac{2^\sfrac{2}{p}}{\sigma_{0}^\sfrac{p+1}{p}} \kap{c} \E{\log(\Sigma_{k}) - \log(\sigma_{0})}  \Biggr).
 \end{align*}
 Using Jensen inequality, the fact that $\Sigma_k$ is a non-decreasing
 sequence, and the definition of $\kap{e}$ and $\kap{f}$ in
 \eqref{ke-def} and \eqref{kf-def} yields the desired result. 
} 
 
We are now ready to give an upper bound on $\E{\Sigma_k}$, a crucial
step in the theory of adaptive regularization methods (see
\cite{Birgin2016} or  \cite{OFFO-ARp} for instance). We will also need
a result on the solutions of a nonlinear equation that combines
logarithmic, linear, and constant terms. The latter is given in
Appendix~B. 
 
\llem{Sigmakbound}{Suppose that AS.1--AS.5 hold. Then for all $k \geq
	0$, we have that 
	\beqn{exprsigmakbound}
	\E{\Sigma_{k}} \leq \sigma_{\max} \eqdef -(p+1)!(\kap{a} \kap{f} + m \kappa_D \chi_p^2  \kap{f})
	\,W_{-1}\left( \frac{-e^{- \kap{g}}}{(p+1)!(\kap{a}\kap{f} +m \kappa_D \chi_p^2  \kap{f})} \right),
	\eeqn 
	where $\kap{a}$ is given by \req{ka-kb-def}, $\kap{f}$ by \req{kf-def} and
	\beqn{kg-def}
	\kap{g} = \frac{\Gamma_0 + \frac{\Sigma_{0}}{(p+1)!} + \kap{a}\kap{e} +m \kappa_D \chi_p^2  \left(\frac{m + 1}{2} + \kap{e}\right)}
	{\kap{a}\kap{f} + m \kappa_D \chi_p^2  \kap{f}}
	\eeqn
	with
	\beqn{Gamma0def}
	\Gamma_0 \eqdef \E{f(X_0) - f_{\rm low}}
	\eeqn
	and $\kap{e}$ given by \req{kf-def}.
}

\proof{Let $j \in \iiz{k}$.	First see that from \eqref{sigmajnegdef} and the definition of $\xi_j$, we have that
		\begin{align*}
			\sum_{j=0}^{k} \E{\xi_j} = \sum_{j=0}^{k} \sum_{i=1}^{m} \E{\|S_{j-i}\|^{p+1}} &= \sum_{i=1}^{m} \sum_{j=0}^{i-1} \E{\|S_{j-i}\|^{p+1}} + \sum_{i=1}^{m} \sum_{j=i}^{k} \E{\|S_{j-i}\|^{p+1}} \\
			&\leq \sum_{i=1}^{m} i + \sum_{i=1}^{m} \sum_{j=0}^{k-1} \E{\|S_{j}\|^{p+1}} \\
			&\leq \frac{m^2+m}{2} + m \sum_{j=0}^{k-1} \E{\|S_{j}\|^{p+1}}
		\end{align*}
	
	Summing \eqref{Lip-f} for all $j$, taking the
	full expectation and using the tower property, and using the previous inequality to bound $\sum_{j=0}^{k} \E{\xi_j}$, we derive that 
	\begin{align*}
		\E{\sum_{j= 0}^k \frac{\Sigma_{j}}{(p+1)!} \|S_j\|^{p+1}} 
		\leq & \sum_{j=0}^{k} \E{f(X_{j}) - f(X_{j+1})} + \sum_{j= 0}^{k} \kap{a}\E{\|S_j\|^{p+1}} \\ 
		&+  \kappa_D \chi_p^2 \left(\frac{m^2+m}{2} +m \sum_{j=0}^{k-1} \E{\|S_j\|^{p+1}}\right).
	\end{align*}
	
	Using \eqref{vkupdate} to simplify the left-hand side, AS.2, and using \eqref{Esumstep} to bound both
	$\sum_{j=0}^k \E{\|S_j\|^{p+1}}$ and $\sum_{j=0}^{k-1}
	\E{\|S_j\|^{p+1}}$, we obtain that 
	\begin{align*}
		\E{\frac{\Sigma_{k+1} - \Sigma_{0}}{(p+1)!}} 
		\leq & \E{f(X_0)} - f_{\rm low} +  \kap{a} (\kap{e} + \kap{f} \log(\E{\Sigma_{k+1}}) \\ 
		&+ m \kappa_D \chi_p^2 \left(\frac{m + 1}{2}+\kap{e} +\kap{f} \log(\E{\Sigma_{k}})\right).
	\end{align*}
	Using now the definition of $\Gamma_0$ in \eqref{Gamma0def},
	the fact that the $\Sigma_j$ sequence is increasing, the last
	inequality gives that
	\begin{align}\label{Sigmakineq}
		\E{\frac{\Sigma_{k+1}}{(p+1)!}} \leq & \Gamma_0 + \frac{\Sigma_{0}}{(p+1)!} + \kap{a} (\kap{e} + \kap{f} \log(\E{\Sigma_{k+1}})) \nonumber \\ 
		&+ m \kappa_D \chi_p^2 \left(\frac{m + 1}{2}+\kap{e} +\kap{f} \log(\E{\Sigma_{k+1}})\right).
	\end{align}  
	Now define
	\beqn{ggu}
	\gamma_1 \eqdef \kap{a} \kap{f} +m \kappa_D \chi_p^2  \kap{f},
	\ms \gamma_2 \eqdef -\frac{1}{(p+1)!},
	\ms u \eqdef \E{\Sigma_{k+1}},
	\eeqn
	\beqn{ggu1}
	\gamma_3 \eqdef \Gamma_0 + \frac{\Sigma_{0}}{(p+1)!} + \kap{a}\kap{e} + m \kappa_D \chi_p^2\left(\frac{m + 1}{2} + \kap{e}\right) 
	\eeqn
	and observe that that $-3 \gamma_2 <  \gamma_1$ since
	$(p+1)! \kap{a} \geq L_p \geq 3$ and $\kap{f} \geq 1$.
	Define the function \\
	$\psi(t) \eqdef \gamma_1 \log(t) + \gamma_2 t + \gamma_3$.
	The inequality \eqref{Sigmakineq} can then be rewritten as
	\beqn{psirewrite}
	0 \leq \psi(u).
	\eeqn
	The constants $\gamma_1$, $\gamma_2$ and $\gamma_3$ satisfy the
	requirements of Lemma~\ref{solloglinconst} and $\psi$
	therefore admits two roots $\{u_1,u_2\}$ whose expressions are given
	in \eqref{Lambertsol}. Moreover, \eqref{psirewrite} is valid only for
	$u \in [u_1, \, u_2]$. Therefore, we obtain from \eqref{ggu},
	\eqref{ggu1} and \eqref{Lambertsol} that 
	\[
	\E{\Sigma_{k+1}} \leq -(p+1)!\gamma_1 W_{-1}\left( \frac{-1}{(p+1)!\gamma_1} e^{-\frac{\gamma_3}{\gamma_1}}\right).
	\]
	We then derive the desired result because the last
	inequality holds for all $k \geq 0$ and $\Sigma_{k}$ is increasing. 
} 

We now discuss the bound obtained \eqref{exprsigmakbound}. First note
that  it is possible to give a more explicit bound on $\sigma_{\max}$
by finding an upper bound on the value of the involved Lambert
function. This can be obtained by using \cite[Theorem~1]{Chat13} which
states that, for $x>0$, 
\beqn{Lamb-bound}
\left|W_{-1}(-e^{-x-1})\right| \leq 1+ \sqrt{2x} + x.
\eeqn
Remembering that, for $\gamma_1$ and $\gamma_2$ given by \req{ggu},
$\log \left( (p+1)! (\kap{a}\kap{f} + m \kappa_D \chi_p^2  \kap{f}) \right) \geq \log(3) > 1$
and taking $x = \kap{g}  - 1  + \log \left( (p+1)! (\kap{a}\kap{f} + m \kappa_D \chi_p^2  \kap{f}) \right)$
in \req{Lamb-bound} then gives that
\begin{align*}
	\left|W_{-1}\left( \frac{-1}{(p+1)!(\kap{a} \kap{f} + m \kappa_D \chi_p^2  \kap{f})} e^{- \kap{g}}\right)\right|
	&\leq \kap{g} + \log \left((p+1)!(\kap{a}\kap{f} + m \kappa_D \chi_p^2  \kap{f} )\right)\\  
	&+\sqrt{2\left( \kap{g} + \log \left( (p+1)! (\kap{a} \kap{f} + m \kappa_D \chi_p^2  \kap{f}) \right) - 1\right)}.
\end{align*}

Complexity results for adaptive regularization methods typically bound
the norm of the gradient at a specific iteration (see
\cite{Birgin2016} for instance), but our differs in this respect and instead
uses the structure of \eqref{vkupdate}, property \eqref{crucial2}, and
Lemma~\ref{Sigmakbound} to produce a bound involving the history of
past gradients.

\lthm{complexity}{
  Suppose that AS.1--AS.5 hold. Then 
  \beqn{Gminzerokbound}
  \min_{j \in \iiz{k}} \E{\|G_{j+1} \|}
  \leq (\kap{c} + \kap{d}m)^\sfrac{p}{p+1}\frac{\sigma_{\max}}{(k+1)^\sfrac{p}{p+1}}
  \eeqn
  where $\kap{c}$ and $\kap{d}$ are given by \req{kc-kd-def} and $\sigma_{\max}$ is
  given by \req{exprsigmakbound}.
}

\proof{
Let $k \geq 1$ and $j \in \iiz{k}$. Taking the full expectation in \eqref{crucial2} and using the tower property, simplifying the upper-bound, using that $\sum_{j= 0}^k \sum_{\ell = j-m}^{j-1} \Sigma_{j+1}-\Sigma_{j} \leq m \Sigma_{k}$ from \eqref{vkupdate}, and using Lemma~\ref{Sigmakbound}, we derive that
\begin{align}\label{sumGjbound}
\bigsum_{j=0}^k \E{\frac{\|{G_{j+1}}\|^\sfrac{p+1}{p}}{\Sigma_{j+1}^\sfrac{1}{p}}  }
\leq\ &  \kap{c} \bigsum_{j=0}^k  \E{\Sigma_{j+1}-\Sigma_{j}} 
+ \kap{d}  \bigsum_{j=0}^k \sum_{\ell = j-m}^{j-1} \E{\Sigma_{\ell+1}-\Sigma_{\ell}} \nonumber \\
&\leq \kap{c} \E{\Sigma_{k+1}} + \kap{d}m \E{\Sigma_{k}} \nonumber \\*[1.5ex]
&\leq  \left( \kap{c} +\kap{d}m \right) \sigma_{\max}.
\end{align}

We now derive a lower-bound on the left-hand side of the last inequality. 
From the H\"older inequality with $q = \frac{p+1}{p}$ and $r = p+1$
and the fact that \eqref{exprsigmakbound} holds, we obtain that
\begin{align}\label{excepglone}
\E{\|G_{j + 1} \|}
&= \E{\frac{\|G_{j + 1} \|}{\Sigma_{j+1}^\sfrac{1}{p+1}}\Sigma_{j+1}^\sfrac{1}{p+1} }
\leq \left(\E{\frac{\|G_{j + 1} \|^\sfrac{p+1}{p} }{\Sigma_{j+1}^\sfrac{1}{p}} }\right)^\sfrac{p}{p+1} \E{\Sigma_{j+1}}^\sfrac{1}{p+1} \nonumber \\ 
&\leq \left(\E{\frac{\|G_{j + 1} \|^\sfrac{p+1}{p} }{\Sigma_{j+1}^\sfrac{1}{p}} }\right)^\sfrac{p}{p+1} \sigma_{\max}^\sfrac{1}{p+1}. 
\end{align}
Taking the last inequality to the power $\sfrac{p+1}{p}$ and using the
result to find a lower bound on the left-hand side of
\eqref{sumGjbound} yields that
\begin{align*}
\frac{\min_{j \in \iiz{k}} \E{\|G_{j+1} \|}^\sfrac{p+1}{p}(k+1)}{\sigma_{\max}^\sfrac{1}{p}}
&\leq \sum_{j=0}^k \frac{\E{\|G_{j + 1} \|}^\sfrac{p+1}{p}}{\sigma_{\max}^\sfrac{1}{p}} 
\leq  \left(\kap{b} + \kap{d} m  \right) \sigma_{\max}.
\end{align*}

Rearranging this last inequality and taking the
$(\sfrac{p+1}{p}$)-th root finally gives the desired result.
} 

The order of dependence on $\epsilon$ given by
Theorem~\ref{complexity} is consistent with that presented in 
\cite{Birgin2016} for the deterministic adaptive regularization algorithm
\cite{Birgin2016,CartGoulToin22}, which has been shown to be optimal for $p$th order
nonconvex optimization \cite{Carmon2019a}. It is also consistent, from
this point of view, with that proposed in \cite{OFFO-ARp} for the deterministic version
of our algorithm.
Theorem~\ref{complexity} however slightly improve on this latter result in another
respect: because the present paper uses different and sharper bounding
techniques, the dependence of $\sigma_{\max}$ on $L_p$ in the constants of \req{Gminzerokbound} is
now $\calO(L_p^{(2p+1)/p}\log(L_p))$, while
that stated in \cite{OFFO-ARp} is $\calO(L_p^{(3p+1)/p})$. 

While the last theorem covers all model degrees, it is worthwhile to
isolate the cases where $p$ is either 1 or 2, detailing some of the 
constants hidden in \req{Gminzerokbound}. We start with
$p=1$.

\lcor{OFFOpone}{Suppose that AS.1--AS.3 and AS.5 hold and that
  $p=1$. Then, the gradients of the iterates generated by
  Algorithm~\ref{InexOFFARp} verify 
\[
\min_{j \in \iiz{k}} \E{\|G_{j+1} \|} \leq \sqrt{ \left(4 L_1^2 + 2\theta_1^2 \sigma_{0}^2 + 2^{m+2} \kappa_D m \right)} \frac{\sigma_{\max}}{\sigma_{0}\sqrt{(k+1)}}
\]
where $\sigma_{\max}$ is defined in \eqref{exprsigmakbound}.
}

Thus obtaining an iterate satisfying $\E{\|G_{k+1} \|} \leq \epsilon$,
 requires at most $\mathcal{O}\left( \epsilon^{-2}\right)$ iterations,
 achieving the complexity rate of linesearch steepest descent
 \cite{Cartis2022-od}. This result is not surprising, since our
 condition \eqref{tensderiverror} for $p=1$ is very similar to the
 strong growth condition \cite{SchLeroux13}. Well-tuned stochastic
 gradient descent reaches the complexity
 rate of deterministic first-order methods under this condition. See
 \cite{khaled2022better} for more details on the theory of stochastic
 gradient descent for nonconvex functions. 

For $p=2$, Theorem~\ref{complexity} may be rephrased as follows.

\lcor{OFFOptwo}{Suppose that AS.1--AS.5 hold and that $p=2$. Then the
  gradients of the iterates generated by Algorithm~\ref{InexOFFARp} verify
	\[
	\min_{j \in \iiz{k}} \E{\|G_{j+1} \|} \leq \sqrt[3]{2 \left( \frac{L_2^\sfrac{3}{2}}{\sqrt{2}} + \frac{\sqrt{2}}{2} + \frac{\theta_1^\sfrac{3}{2}}{2^\sfrac{3}{2}} + 2^{m-1}(4+\sqrt{2}) \kappa_D m \right)^2}  \frac{\sigma_{\max}}{\sigma_{0}(k+1)^\sfrac{2}{3}}
	\]
	where $\sigma_{\max}$ is defined in \eqref{exprsigmakbound}.
}

Again, if we are interested in reaching an iterate such that
$\E{\|G_{k+1}\|} \leq \epsilon$,
$\mathcal{O}\left(\epsilon^{-3/2}\right)$ iterations are required in
the worst case, achieving the same rate as optimal second-order
methods (see \cite{Cartis2022-od} and the references therein). As a
consequence, our algorithm is an optimal adaptive cubic regularization
method without function evaluation in a fully stochastic setting.

\numsection{Applications of the StOFFAR$p$ algorithm}\label{applicationsstooffo}

In this section, we present a series of practical cases in support of
the StOFFAR$p$ algorithm. Given that the analysis of the latter is
contingent on AS.5, two frameworks are provided that satisfy this
condition.  We start by considering inexact derivatives in the context
of multiprecision arithmetic. We then provide conditions on subsample
size of the stochastic gradient and Hessian for practical machine
learning problems with $p=2$.

\subsection{Inexact Derivatives}
  
\subsection{Inexact Derivatives}

Our theory naturally applies to the case where derivatives are
inexact. For the sake of clarity, we drop the uppercase notation and
use only lowercase in this subsection. For this particular case,
\eqref{tensderiverror} holds without expectation for all
iterations. Specifically, there exists $\kappa_D > 0$ such that the
inaccurate derivatives $\overline{\nabla_x^i f}(x_k)$ used to compute
the model \eqref{model} satisfy,  
\beqn{tensderivEDA}
\| \nabla_x^i f(x_k) - \overline{\nabla_x^i f}(x_k) \|
\leq \kappa_D  \sum_{j=1}^{m}\|s_{k-j} \|^{p+1-i} \tim{ for all } i \in \ii{p}.
\eeqn

The conditions here are very similar to those proposed in
\eqref{tensderivcurrstep}. Again, one of the advantages of
\eqref{tensderivEDA} is that it considers the previous steps and not
the current one, allowing \eqref{tensderivEDA} to be enforced at the
beginning of each iteration. This approach formally covers the use of
imprecise derivatives, where the approximation of high-order tensors
is performed by using finite differences of low-order
derivatives. For more details on these algorithmic variants, we refer
the reader to \cite[Subsection 13.2]{Cartis2022-od}. 

The inexact version of our algorithm also falls under the Explicit
Dynamic Accuracy (EDA) framework \cite[Section~13.3]{CartGoulToin22},
since the conditions can be enforced a priori. The aforementioned
settings are a hot topic, and algorithms have recently been proposed
\cite{BellavGraRicc18,GrattonToin20}. These theoretical advances have
arisen to take advantage of developments in large-scale modern
computing hardware that allow loose numerical approximations of
derivatives when needed.  An imprecise version of our \al{StOFFAR$p$}
can be used in this context and may even offer a simpler alternative
compared to current explicit dynamic accuracy adaptive regularization
methods (see for example \cite[Algorithm 13.3.3]{Cartis2022-od}).
 
\subsection{Machine Learning Problems}\label{ML-probs}

In this subsection, we focus on the case where $p=2$ as the results of
this section are focused on practical machine learning problems. In
the latter case \eqref{problem} becomes  
\beqn{probML}
\min_{x \in \Re^n} \{ f(x) = \frac{1}{N} \sum_{i=1}^{N}  f_i(x,y_i,a_i) \}
\eeqn
where both $N$ and $n$ (the number of optimized variables) may
exceed to millions and $f_i$ may be nonconvex. The pairs  $(a_j,y_j)$
are independent and identically distributed random variables coming
from an a priori unknown distribution $\mathcal{D}$.  In this case, it
is common to randomly sample  batches of indices in the expression of
$f$ to approximate its derivatives. The sampled gradient and Hessian
are therefore given by
\beqn{gradHesssampling}
 \overder{1}{k} = \frac{1}{b_{g,k}} \sum_{i \in \calB_{g,k}}\nabla_x^1 f_i(X_k,y_i,a_i),
 \tim{ and }
 \overder{2}{k} = \frac{1}{b_{H,k}} \sum_{i \in \calB_{H,k}} \nabla_x^2 f_i(X_k,y_i,a_i),
\eeqn
where $\calB_{g,k}$ and $\calB_{H,k}$ are the batches at iteration $k$
of cardinality $b_{g,k}$ and $b_{H,k}$, respectively. In our case,
$\nabla_x^1 f_i(X_k,y_i,a_i)$ are i.i.d.\footnote{independent
identically distributed} vector-valued random variables and
$\nabla_x^2 f_i(X_k,y_i,a_i)$ are i.i.d. random self-adjoint matrices
with dimension $n \times n$.  To obtain lower bounds on batch sizes
$b_{g,k}$ and $b_{H,k}$ of the stochastic gradient and Hessians
\eqref{gradHesssampling}, conditions in expectation on the noise of
the gradient and Hessian of each $f_i$ must be assumed. For clarity,
we will drop the $\Ek{.}{k}$ notation and keep only $\E{.}$ since we
focus only on a specific iteration $k$.  The goal of the next theorem
is to provide requirements on $b_{g,k}$ and $b_{H,k}$ under
assumptions that are common in the literature
\cite{ChaytiDoiJaggi23,ZhouXuGu19} in order to satisfy
\eqref{tensderiverror}.

\lthm{samplingcond}{ Let $k$ be an iteration of the \al{StOFFAR$2$}
  algorithm and suppose that the objective function has the structure
  given in \eqref{probML} and that for each $i \in \ii{\bigNred}$, there
  exist non-negative constants $\sigma_{g}$ and $\sigma_H$ such that 
\beqn{excpeHessgraderror}
\E{\| \nabla_x^1f_i(X_k,y_i,a_i)-\nabla_x^1f(X_k)\|^2} \leq\sigma_g^2
\tim{ and }
\E{\| \nabla_x^2f_i(X_k,y_i,a_i)-\nabla_x^2f(X_k)\|^3} \leq \sigma_H^3.
\eeqn
Then the estimators introduced in \eqref{gradHesssampling} for problem
\eqref{probML} verify conditions \eqref{tensderiverror} if  
\beqn{bgksamplebound}
b_{g,k} \geq \frac{\sigma_g^2}{\kappa_D^\sfrac{4}{3} \xi_k^\sfrac{4}{3}},
\eeqn 
and 
\beqn{bHksamplebound}
b_{H,k} \geq \frac{9 \sigma_H^2 e \log(n)}{2\kappa_D^\sfrac{2}{3} \xi_k^\sfrac{2}{3}}. 
\eeqn
}

\proof{As the proof combines elements already developed in
  \cite{ChaytiDoiJaggi23,ZhouXuGu19} but adapted to take into account
  \eqref{tensderiverror}, it is differed to
  Appendix~\ref{samplingproof}. } 

Before proceeding, we discuss our proposed sampling conditions and
provide a discussion when  $ m = 1$ and $\xi_k =
\|S_{k-1}\|^3$. First, note that we have obtained the same order of
dependence on the step size as in the work of \cite{kohler17a}. Our
framework improves on this reference because we have not imposed
Lipschitz continuity on $f_i$ or its derivatives. Moreover, our
condition covers the use of the previous step to scale the
batch-sizes, whereas the theoretical result developed in
\cite{kohler17a} uses the current step size.
Finally, it should be noted that our framework is more flexible than
previous works \cite{kohler17a,BellaviaGurioli21} in that it allows
the error to depend on the past $m$ steps, rather than just a specific
one. 

\numsection{Numerical illustration}\label{numerics}

In this section, we illustrate the numerical behaviour of our proposed
\al{StOFFAR$p$} algorithm for $p=1$ and $p=2$ for the machine-learning
problems discussed in Subsection~\ref{ML-probs}. The 
goal of the following experiments is to demonstrate the advantages of
high-order objective-free function algorithms for machine-learning
problems. We perform numerical tests on two different formulations of
the binary classification problem. Throughout this section, $\{a_i, \,
y_i\}_{i=1}^{\bigNred}$ represents the training data with $a_i \in \Re^n$ and
$y_i \in \{0,\,1\}$ representing the $i$th feature and the $i$th
target label, respectively. For the binary classification, we propose
the following formulation as a minimization task: 
\beqn{sigmoidbin}
\min_{x \in \Re^n} f(x)
= \min_{x \in \Re^n} \frac{1}{N} \sum_{i=1}^{N}  f_i(x,y_i,a_i)
= \min_{x \in \Re^n} \frac{1}{N} \sum_{i=1}^{N} \left(y_i - \phi(a_i^\intercal x) \right)^2,
\eeqn
where 
\beqn{phisigmoid}
\phi(a^\intercal x) = \frac{1}{1+e^{-a^\intercal x}}.
\eeqn
This minimization problem has already been considered in
\cite{BellaviaGurioli21,BellaviaGianBen2020}.
We refer the reader to these references for the expressions of both
the gradient and the Hessian. We also consider a second case of
nonconvex binary classification studied in \cite{kohler17a,
  khaled2022better}, where a standard binary 
logistic regression is regularized with a nonconvex term. The binary
classification problem is then formalized as: 
\beqn{binarylogisticreg}
\min_{x \in \Re^n} f(x) = \min_{x\in \Re^n} \frac{1}{N} \sum_{i=1}^{N} \left( -y_i \log(\phi(a_i^\intercal x)) - (1-y_i) \log(1-\phi(a_i^\intercal x)) + \alpha \sum_{j=1}^{n}
 \frac{x_j^2}{1+x_j^2} \right),
\eeqn 
where $\alpha$ is a parameter that regulates the strength of the
penalization. Note that from \cite[Subsection~2.1]{Drusvyatskiy18} and by using a judicious decomposition of both  functions \eqref{sigmoidbin} or \eqref{binarylogisticreg}, any subsubsampled Hessian of \eqref{sigmoidbin} or \eqref{binarylogisticreg} satisfies Assumption~\ref{assumption:AS4}. The rest of this section is organized as follows.
Implementation issues are considered in
Subsection~\ref{implementation}. To satisfy AS.5, we run a variant of
our \al{StOFFAR$p$} with $p=2$ and various $m$ values, denoted
\texttt{OFFAR$2$-m}, that implements a sampling strategy using the
scaling rules given in Theorem~\ref{samplingcond}. As a baseline, we
use a \al{StOFFAR$p$} with $p=1$ and $m=1$, denoted \texttt{WNGRAD},
since our algorithm retrieves the method proposed in
\cite{WNGRAD}. As in Theorem~\ref{samplingcond}, we also derive a
condition on the sample size to satisfy AS.5 for this method. We have
avoided comparison with other second-order stochastic algorithms, such
as those proposed in \cite{BellaviaGianBen2020, BellaviaGurioli21,
  kohler17a, WangZhouYingbinLan19}, since they either require access
to the exact value of the function to adjust the regularization
parameter $\sigma_k$, or assume knowledge of the Lipschitz
constant. Some illustrations of both methods are provided in
Subsection~\ref{res-subs}. 

\subsection{Implementation Issues}\label{implementation}

Our implementation relies on the code provided in
\cite{kohler17a}\footnote{Available at
{https://github.com/dalab/subsampled\_cubic\_regularization}.}. The
subsampled cubic regularization subroutine is slightly adapted to
allow the use of the update rule given in \eqref{vkupdate}, to fulfill
the condition given in \eqref{gradstep} when computing the step, and
to subsample in accordance with the conditions of
Theorem~\ref{samplingcond}. Specifically, at the initial iteration of
the \texttt{OFFAR$2$-m} algorithm, the values of $b_{h,0}$ and
$b_{g,0}$ are set to $0.05 \cdot M$ and $0.20 \cdot N$ in order to
compute the approximate Hessian and gradient, as defined in
\eqref{gradHesssampling}. Note that this choice of initial subsampling
size is consistent with past subsampled methods developed in the
literature \cite{kohler17a, BellaviaGurioli21}. For $k \geq 1$, we
using the following subsampling strategy: 
\beqn{subsamplinkiter}
b_{g,k} = \max\left (\frac{c_g}{\xi_k^\sfrac{4}{3}}, 0.20\cdot N \right ), \quad 
b_{H,k} = \max \left (\frac{c_H}{\xi_k^\sfrac{2}{3}}, 0.05 \cdot N \right ),
\eeqn 
where $c_g = b_{g,0}m^\sfrac{4}{3}$ and $c_H = \frac{b_{h,0}
  m^\sfrac{2}{3}}{\log(n)}$ and $m$ is defined in
\eqref{sigmajnegdef}. The choices of the constants $c_g$ and
$c_H$ are made to ensure that our first subsampled derivatives verify
\eqref{bgksamplebound} and \eqref{bHksamplebound} with $\kappa_D$
chosen as $\frac{\sigma_{g}^\sfrac{3}{2}}{b_{g,0}^\sfrac{3}{4}m}$ for
gradient subsampling and
$\frac{(9e)^\sfrac{3}{2}(\log(n)\sigma_H)^\sfrac{1}{3}}{(2b_{H,0})^\sfrac{3}{2}m}$
for the Hessian subsampling. We also chose $\sigma_{0} = 0.01$ and
$\theta_1 = 2$ and ran four variants with $\texttt{m} \in
\{1,50,250,500\}$. 

We also developed our own implementation of the \texttt{WNGRAD}
algorithm where we use an initial batch size of $b_{g,0} = 0.05\cdot
N$, $m = 1$, and subsample for $k \geq 1$ with  
\beqn{WNGRADsubsamples}
b_{g,k} = \max \left(0.05 \cdot N, \frac{0.1}{\|S_{k-1}\|^2}\right).
\eeqn 
We also choose $\sigma_0 = 0.1$ for \texttt{WNGRAD}. Both methods
start from an initial point $x_0 = (0, \, 0, \, \dots, 0)$ and
$\alpha$ in \eqref{binarylogisticreg} is taken equal to $0.001$. 

The algorithms are stopped when an iterate $x_k$ satisfying
\beqn{stopcriteria}
\|\overline{\nabla_x^1 f}(x_k) \|
\leq \epsilon \quad \text{with} \quad \epsilon = 0.0005
\eeqn
is reached. The maximum number of iterations for both
\texttt{OFFAR2-m} and \texttt{WNGRAD} is set to $1000$ and $10000$,
respectively. The datasets are taken from the LIBSVM library
\cite{LIBSVM} (see Appendix~\ref{Datasets} for more detail).

\subsection{Results}\label{res-subs}

To evaluate the performance of our methods that involve stochastic
ingredients (resulting from approximation by subsampling), all
reported results are averages over $20$ independent runs. To provide
an appropriate comparison between the tested methods which may employ
different batch sizes, we report the performance measure 
\beqn{complextau}
\tau_{algo} = \sum_{i=1}^{k} (b_{g,i}+b_{H,i}) \cdot \textrm{ege}_i,
\eeqn
where $\textrm{ege}_i$, the effective gradient evaluation metric,
counts the number of Hessian-vector products used at each iteration to
compute the step for the \texttt{OFFAR2-m} methods in addition to the
number of gradient evaluations. For the \texttt{WNGRAD} algorithm, the
value of $\textrm{ege}_i$ is equal to one, while the value of
$b_{H,i}$ is equal to zero.  

\begin{figure}[!ht]
	\centering
	\includegraphics[width=10cm, keepaspectratio]{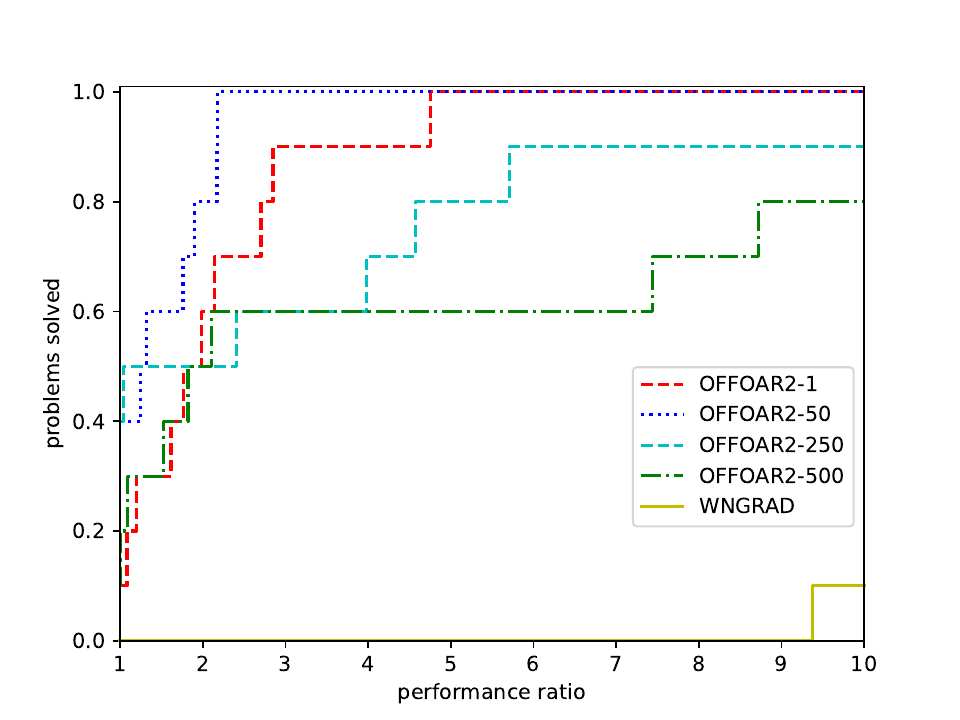}
	\caption{\label{perfprof} Performance profile of \texttt{OFFAR2-m}
      for {\tt SUSY} and {\tt w8a} for $\texttt{m} \in \{1,50,250,500\}$ and \texttt{WNGRAD}}
\end{figure}

\begin{figure}[!ht]
	\centering
	\includegraphics[width=7cm]{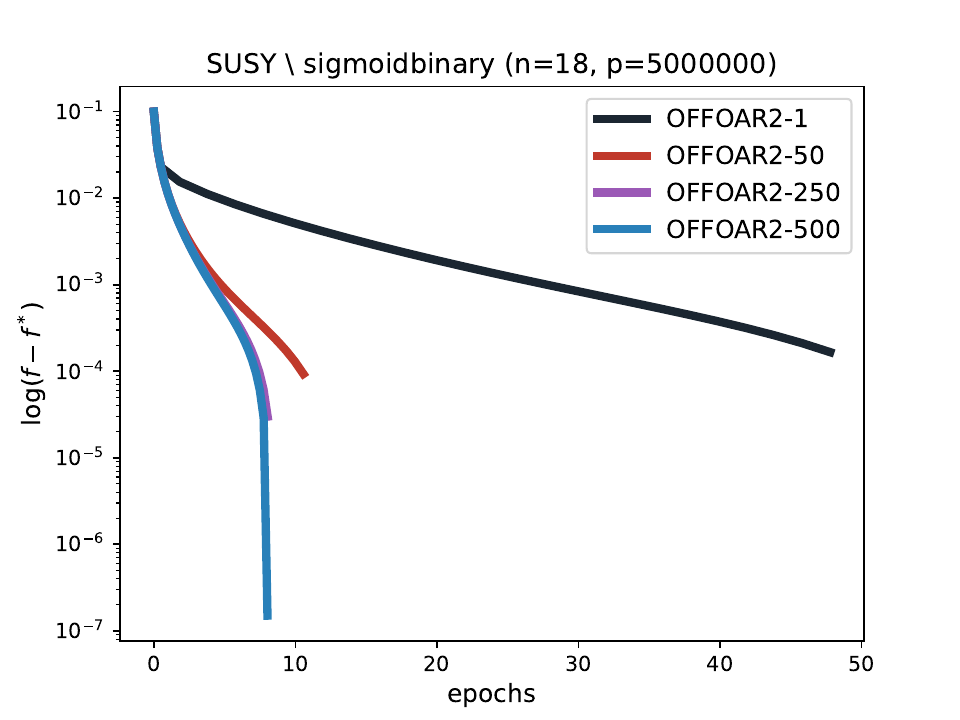}
	\hspace*{-5mm}
	\includegraphics[width=7cm]{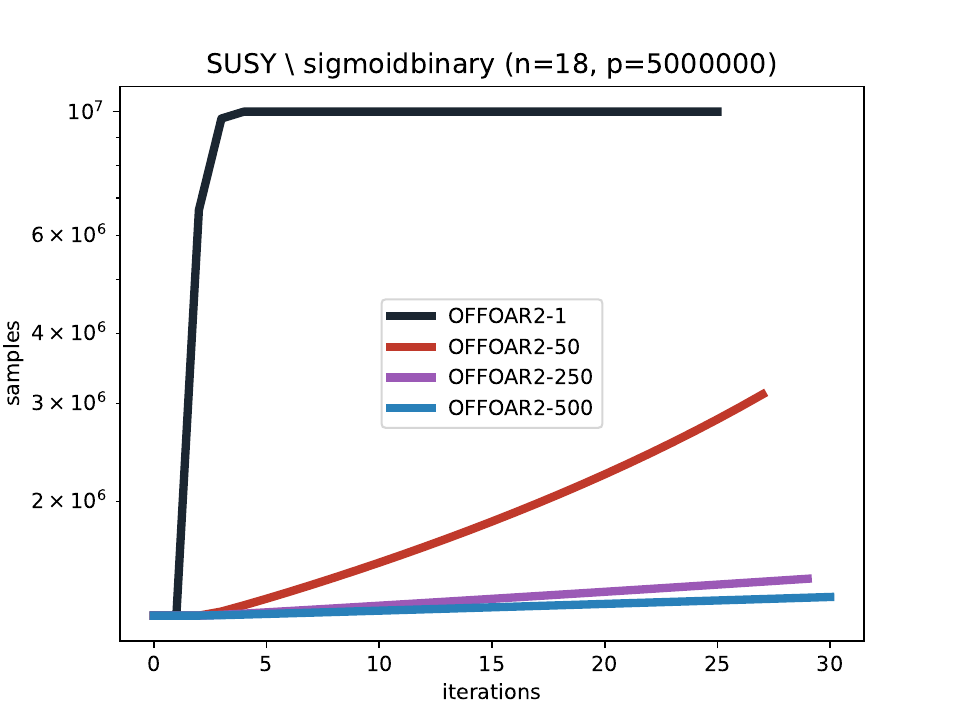}
	\hspace*{-5mm}
	\includegraphics[width=7cm]{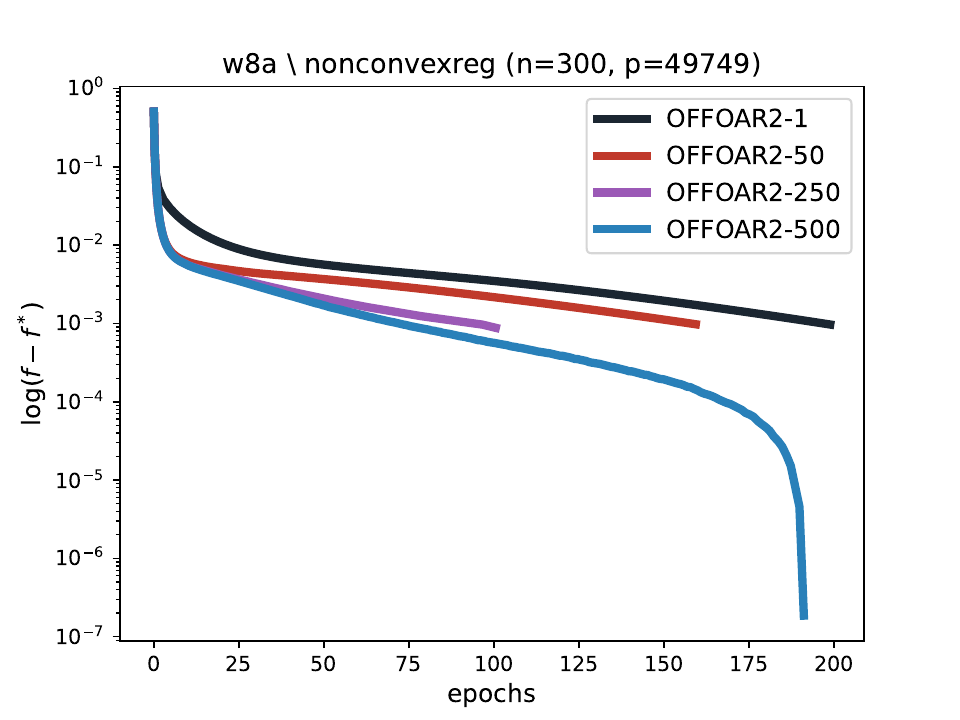}
	\hspace*{-5mm}
	\includegraphics[width=7cm]{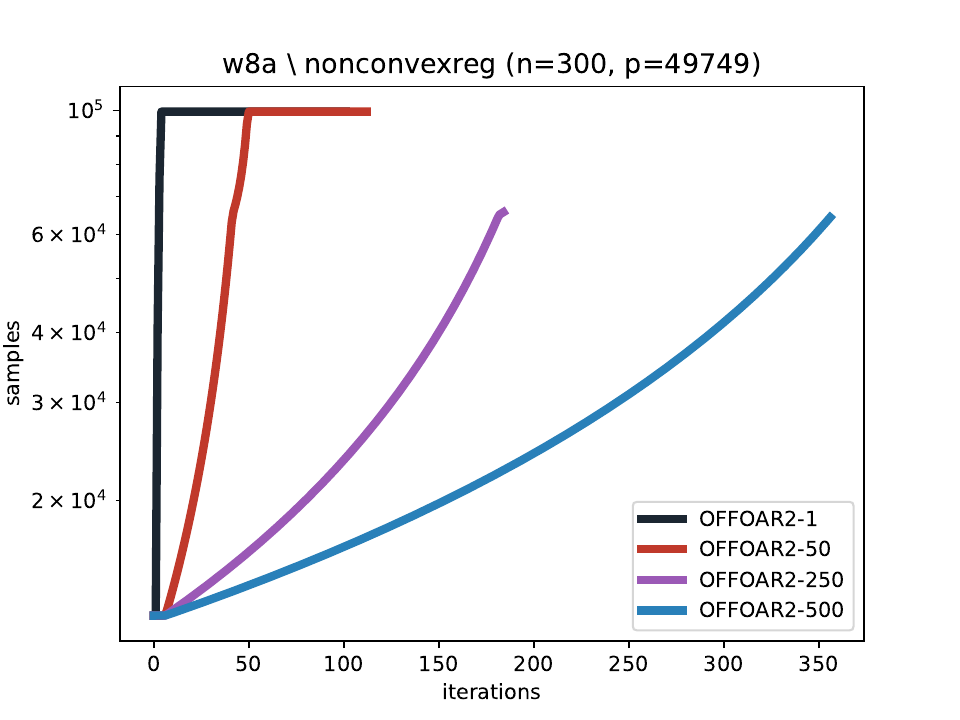}
	\hspace*{-5mm}
	\caption{\label{fig:objf} Evolution of the loss function w.r.t the
      epochs and number of samples along iterations for {\tt SUSY} and
      {\tt w8a}}
\end{figure}

Figure~\ref{perfprof} shows the standard performance profile
\cite{DolaMoreMuns06} for the five methods with respect to the
performance measure \eqref{complextau}. The figure
illustrates that our proposed second-order method identifies an
approximate first-order stationary point more rapidly than simple
adaptive gradient methods represented here by \texttt{WNGRAD}. It is
evident that the \texttt{OFFAR2-50} method is the most
efficient. \texttt{OFFAR2-1} is the second most efficient, but it may
perform less effectively than the other methods on some problems.  To
illustrate this point, consider Figure~\ref{fig:objf}, in which
$f^\star$ denotes the best (i.e., minimum) value obtained among all
the four tested methods, and where the number of samples is reported
for each method.  The two problems considered here, \texttt{SUSY} and
\texttt{w8A}, illustrate cases where the performance of the methods
with longer memory ($\texttt{m} \in \{250, 500\}$) is superior to that
of the methods with shorter memory ($\texttt{m} \in \{1, 50\}$), both
in terms of convergence and the number of samples used by the methods.

We remind the reader that one epoch denotes a pass made on the whole
data set (samples $f_i$ in \eqref{probML}) when computing the
stochastic gradient and the Hessian.

We also use Figure~\ref{fig:objf} to exemplify a generic problem
occuring when using short memory and single-step length control: the
resulting method may become exact (and therefore computationally
expensive) after only a few iterations, as the required samples
involve the entire data set.  This (undesirable) behavior is also
observed in many methods, including the subsampled cubic
regularization method presented in \cite{kohler17a}. Other subsampled
second-order methods impose a tight probabilistic bound on all
iterations (as shown in \cite{RoostaKhorasani18, Yao2021,
  ChenJianetal22}), which again causes the algorithm to be
deterministic for most iterations. The same drawback also appears in
some methods that impose the growth of the batch size and the use of
the exact Hessian and gradient starting from a specific iteration
\cite{Bollapragadaetal18}.  In contrast, \texttt{OFFAR2-m} methods
with long memory allow the error bounds to remain large, resulting in
more aggressive sampling until termination. It is worth noting that,
empirically, \texttt{OFFAR2-m} with large \texttt{m} reaches local
minima with a lower objective value.  Longer memory may however
occasionally may result in slower convergence in practice (as
illustrated by Figure~\ref{perfprof}), and satisfying the criteria
\eqref{stopcriteria} may become costly. The reader is referred to the
examples shown in Appendix~\ref{addfig} to understand some of the
problems that arise when using a large \texttt{m}.

These early results suggest that using high-order OFFO algorithms may
be beneficial, but the authors are aware that additional numerical
experiments are required to better assess their potential. Indeed,
refinements on the update rule of the regularization parameter have
been proposed in OFFO second-order methods, be it trust-region
\cite{Grapiglia2022} or adaptive regularization \cite{OFFO-ARp}, and a
thorough analysis of the influence of the values of $\|S_{-1}\|,
\ldots, \|S_{-m}\|$ and the length of the memory \texttt{m} may be
required. We have avoided their discussion here to keep our analysis
and numerical experiments concise. Their addition may require a more
involved proof and may impose stronger assumptions. 

\numsection{Discussion}\label{dicuss-s}

In this paper, we have developed a fully stochastic theory for an
objective-function-free adaptive regularization algorithm described in
\cite{OFFO-ARp}. Since the algorithm does not use the function value
to accept or reject the step, it avoids the need to compute this value
with an accuracy higher than that used for the gradient, thereby
making it a computationally attractive technique for noisy problems.
The new algorithm introduces novel conditions on the probabilistic tensor
derivatives, and uses the history of past steps to determine the level
of derivatives' accuracy which is acceptable in expectation to ensure
convergence. Our analysis shows that its evaluation complexity is
optimal in order.

We also discussed two application cases. The first focuses on noisy
inexact functions, where inaccuracy arises from lower precision
computations or the use of finite differences. The second case is
finite-sum minimization (typical of machine learning problems), where
we provide sample size conditions to meet the specified requirements
under mild assumptions. Applying the algorithm to practical binary
classification problems highlighted the advantage of second-order OFFO
methods over standard adaptive gradient strategies and also showed that
the proposed sampling scheme can remain practical throughout the
computation.

Unsurprisingly, an extension of the algorithm to guarantee termination
at approximate second-order stationary points is possible, in the vein
of what was proposed in \cite[Section 4]{OFFO-ARp} for the
deterministic case. The analysis would be very similar to that of
Section~\ref{complexity-s}, replacing $\|G_k\|$ by the appropriate
measure of criticality.

One possible further improvement is to study the OFFO algorithm under
the  assumption that 
\[
\Ek{\| \nabla_x^i f(X_k) - \overline{\nabla_x^i f}(X_k) \|^\sfrac{p+1}{p+1-i}}{k} \leq \kappa_D \|S_{k-1}\|^{p+1} + \kappa_c, \quad \text{for all } i \in \{1, \ldots, p\}.
\]
An assumption of this nature has been considered in the analysis of adaptive
gradient methods \cite{Fawetal22}, and extending it to higher-order
OFFO schemes seems a natural line for future research.  
A second line may focus on proposing OFFO schemes that incorporate
momentum when updating the regularization parameter.

\appendix
\vspace*{2mm}
\section{Proof of \eqref{skbound} }\label{skborne}
\proof{ In the following, we use lowercase notation as the Lemma~\ref{stepgkbound} is valid for all iterations and all realizations.
	If $p=1$, we obtain from \req{descentmodel} and the Cauchy-Schwartz
	inequality that
	\[
	\half \sigma_k \|s_k\|^2 < -\overline{g}_k^\intercal s_k \leq \|\overline{g}_k\|\,\|s_k\|
	\]
	and \req{skbound} holds with $\eta = 0$.
	Suppose now that $p>1$. \req{descentmodel} gives that
	\begin{align*}
	\frac{\sigma_{k}}{(p+1)!} \|s_k\|^{p+1}
	\leq  -\overline{g}_k^\intercal s_k -  \sum_{i=2}^p \frac{1}{i!} \overline{\nabla_x^i f}(x_k) [ s_k]^i  
	\leq \| \overline{g}_k\| \| s_k\|  +   \sum_{i=2}^p \frac{{\khigh}}{i!} \| s_k\|^i,
	\end{align*}
	where we applied  AS.4 to obtain the last inequality.
	
	\noindent 
	Applying now the Lagrange bound for polynomial roots \cite[Lecture~VI, Lemma~5]{Yap99}
	with $x=\|s_k\|$, $n=p+1$, $a_0 = 0$, $a_1=\| \overline{g_k}\|$, $a_i = {\khigh}/i!$ $i \in \iibe{2}{p}$ and $a_{p+1}
	=\sigma_k/(p+1)!$, we obtain from \req{descentmodel} that the equation
	$\sum_{i=0}^na_ix^i = 0$ admits at least one strictly positive root,
	and we may thus derive that
	\begin{align*}
		\|s_k\|\leq  2 \max\Biggl(&\left( \frac{
          (p+1)!\|\overline{g}_k\|}{\sigma_k}\right)^\sfrac{1}{p},
        \left\{ \left[\frac{{\khigh}(p+1)!}{i! \sigma_k }\right]^\sfrac{1}{p-i+1}\right\}_{i\in \iibe{2}{p}} \Biggr).
	\end{align*}
	Using now the fact that $\sigma_{k} \geq \sigma_{0}$ and the
    definition of $\eta$ in \eqref{kaphigh} yields \eqref{skbound}. 
}

\section{Solutions of the equation $\gamma_1 \log(u) + \gamma_2u + \gamma_3 = 0$ }
\renewcommand{\theequation}{\thesection.\arabic{equation}}
\setcounter{equation}{0}
\llem{solloglinconst}{Let $(\gamma_1,\gamma_2,\gamma_3) \in \Re_{+}^{\star} \times \Re_-^{\star} \times \Re^{+}$ and $\frac{\gamma_2}{\gamma_1} \geq -\frac{1}{3} $. Then the equation 
\beqn{tosolveqgeneric}
\gamma_1 \log(u) + \gamma_2 u + \gamma_3 = 0
\eeqn
admits two solutions $0 < u_1 < u_2$ given by
\begin{equation}\label{Lambertsol}
u_1 = \frac{\gamma_1}{\gamma_2} W_{0}\left( \frac{\gamma_2}{\gamma_1}
e^{\frac{-\gamma_3}{\gamma_1}}\right)
\;\mbox{ and }\;
u_2 = \frac{\gamma_1}{\gamma_2} W_{-1} \left(
\frac{\gamma_2}{\gamma_1} e^{\frac{-\gamma_3}{\gamma_1}}\right),
\end{equation}
where $W_0$ and $W_{-1}$ are the two branches of the Lambert function \cite{Corless1996}.
}
\proof{
Note that since $e^{\frac{-\gamma_3}{\gamma_1}} \leq 1$
and $ -\frac{1}{3} \leq \frac{\gamma_2}{\gamma_1} <0 $, we obtain that
\beqn{Wbiendef}
-\frac{1}{3} \leq \frac{\gamma_2}{\gamma_1} e^{\frac{-\gamma_3}{\gamma_1}} < 0.
\eeqn
Let $u$ be a solution of \eqref{tosolveqgeneric}. Rearranging the
equality \eqref{tosolveqgeneric} and taking the exponential yields
that
\[
u= e^{\frac{-\gamma_3}{\gamma_1} - \frac{\gamma_2}{\gamma_1} u }
\]
and thus that
\[
\frac{\gamma_2}{\gamma_1} u e^{\frac{\gamma_2}{\gamma_1} u }
= \frac{\gamma_2}{\gamma_1} e^{\frac{-\gamma_3}{\gamma_1}}.
\]
Taking $w \eqdef \frac{\gamma_2}{\gamma_2} u  $ and using
\eqref{Wbiendef}, we obtain that the equation 
\[
w e^{w} = \frac{\gamma_2}{\gamma_1} e^{\frac{-\gamma_3}{\gamma_1}}
\]
admits two distinct solutions $w_1$ and $w_2$ given by
\[
w_1 = W_{0}\left( \frac{\gamma_2}{\gamma_1} e^{\frac{-\gamma_3}{\gamma_1}}\right), \, \, \quad w_2 = W_{-1} \left( \frac{\gamma_2}{\gamma_1} e^{\frac{-\gamma_3}{\gamma_1}}\right) \tim{ and } w_2 < w_1 < 0.
\]
The desired result then follows from the facts that
$u = \frac{\gamma_1}{\gamma_2} w$
and that $\frac{\gamma_1}{\gamma_2} < 0$.
}
\section{Proof of Theorem~\ref{samplingcond}}\label{samplingproof}
\setcounter{equation}{0}
Before proving Theorem~\ref{samplingcond}, we need the two following auxiliary lemmas that we state below.

\llem{threehalfzi}{Suppose that $z_1, z_2, \dots, z_N$ are i.i.d vector valued random variables with $\E{z_i} = 0$ and $\E{\|z_i\|^2} < + \infty$. Then
	\[
	\E{\left\| \frac{1}{N} \sum_{i=1}^N z_i\right\|^\sfrac{3}{2}} \leq \frac{1}{N^\sfrac{3}{4}} \left( \E{\|z_i\|^2}\right)^\sfrac{3}{4}.
	\]
 }
\proof{See \cite[Lemma~31]{ZhouXuGu19} for the statement of the lemma and Appendix~C of this reference for its proof.}

\llem{matrixconcentration}{ Suppose that $q \geq 2$, $n \geq 2$, and fix $r \geq \max(q, 2 \log n)$. Consider i.i.d. random
	self-adjoint matrices $Y_1, \dots, Y_N$ with dimension $n \times n$, $\E{Y_i} = 0$. Then
\[
\E{ \left\Vert \sum_{i=1}^N Y_i \right\Vert^q}^\sfrac{1}{q} \leq 2 \sqrt{er} \left\Vert  \left(\sum_{i=1}^N \E{Y_i^2}\right)^\sfrac{1}{2} \right\Vert + 4 er \left( \E{\max_{i} \|Y_i\|^q}\right)^\sfrac{1}{q}.
\]
  }
\proof{As for the previous lemma, see \cite[Lemma~32]{ZhouXuGu19} for the statement of the Lemma and its proof. }

We are now in a position to provide the proof of the statement of Theorem~\ref{samplingcond}. 
\proof{
We start by providing a proof on $b_{g,k}$. First, denote  $
g_{i,k} \eqdef \nabla_x^1 f_i(X_k,y_i,a_i) $ for $i \in \calB_{g,k}$,
so that \eqref{gradHesssampling}  and \eqref{probML} give that
\beqn{giksimpfact}
 \E{g_{i,k}} = \nabla_x^1 f(X_k) \tim{ and } \overder{1}{k} = \frac{1}{b_{g,k}} \sum_{i \in \calB_{g,k}} g_{i,k}.
\eeqn

Applying now Lemma~\ref{threehalfzi} with $z_i = \frac{g_{i,k} - \nabla_x^1 f(X_k)}{b_{g,k}}$ for $i \in \calB_{g,k}$ and using the first part of \eqref{excpeHessgraderror}, we derive that
\[
\E{ \left\Vert \frac{1}{b_{g,k}} \sum_{i \in \calB_{g,k}} g_{i,k} - \nabla_x^1 f(X_k) \right\Vert^\sfrac{3}{2}} \leq \frac{\sigma_{g}^\sfrac{3}{2}}{b_{g,k}^\sfrac{3}{4}},
\]
and so if $b_{g,k}$ is taken as in \eqref{bgksamplebound},  \eqref{tensderiverror} holds for $i=1$ and $p=2$.

Again, as for the gradient, we denote $H_{i,k} \eqdef \nabla_x^1
f_i(X_k,y_i,a_i)$., and thus 
\beqn{Hiksimpfact}
\E{H_{i,k}} = \nabla_x^2 f(X_k) \tim{ and } \overder{2}{k} = \frac{1}{b_{H,k}} \sum_{i \in \calB_{H,k}} H_{i,k}.
\eeqn
Also note that \eqref{excpeHessgraderror} and Jensen's inequality imply that
\beqn{Hiktwomoment}
\E{\|H_{i,k} - \nabla_x^2 f(X_k)\|^2} \leq \left(\E{\|H_{i,k} - \nabla_x^2 f(X_k)\|^3}\right)^\sfrac{3}{2} \leq \sigma_H^2.
\eeqn   
Applying now Lemma~\ref{matrixconcentration} with $q = 3$, $r = 2 \log (n)$, $N = b_{H,k}$ and $Y_i = \frac{ H_{i,k} - \nabla_x^2 f(X_k)}{b_{H,k}}$, we obtain that
\begin{align}\label{twotermsHik}
\E{ \left\Vert \frac{1}{b_{H,k}} \sum_{i \in \calB_{H,k}} H_{i,k} -  \nabla_x^2 f(X_k) \right\Vert^3} 
&\leq \Biggl( 2 \sqrt{2e \log(n)} \left\Vert  \left( \sum_{i \in \calB_{H,k}} \frac{1}{b_{H,k}^2} \E{\left(H_{i,k} -  \nabla_x^2 f(X_k) \right)^2}\right)^\sfrac{1}{2} \nonumber \right\Vert\\ 
&+ \frac{8e\log (n)}{b_{H,k}} \left( \E{\max_{i \in \calB_{H,k}} \|H_{i,k} -  \nabla_x^2 f(X_k)\|^3}\right)^\sfrac{1}{3} \Biggr)^3.
\end{align}

Let us now establish a bound on $\left\Vert  \left( \sum_{i \in \calB_{H,k}} \frac{1}{b_{H,k}^2} \E{\left(H_{i,k} -  \nabla_x^2 f(X_k) \right)^2}\right)^\sfrac{1}{2} \right\Vert$. Successively using the fact that $\|A^\sfrac{1}{2}\| = \|A \|^\sfrac{1}{2}$ for any positive definite matrix $A$, the Jensen's inequality, that $\|B^2\| = \|B\|^2$ for any symmetric matrix $B$,  and \eqref{Hiktwomoment}, we derive that
\begin{align}\label{firsttermHik}
\left\Vert  \left( \sum_{i \in \calB_{H,k}} \frac{1}{b_{H,k}^2} \E{\left(H_{i,k} -  \nabla_x^2 f(X_k) \right)^2}\right)^\sfrac{1}{2} \right\Vert 
&=    \left\Vert \sum_{i \in \calB_{H,k}} \frac{1}{b_{H,k}^2} \E{\left(H_{i,k} -  \nabla_x^2 f(X_k) \right)^2} \right\Vert^\sfrac{1}{2} \nonumber \\
&= \left\Vert \frac{1}{b_{H,k}} \E{\left(H_{i,k} -  \nabla_x^2 f(X_k) \right)^2} \right\Vert^\sfrac{1}{2} \nonumber \\
&\leq \frac{1}{\sqrt{b_{H,k}}} \E{ \left\Vert  \left(H_{i,k} -  \nabla_x^2 f(X_k) \right)^2 \right\Vert}^\sfrac{1}{2} \nonumber \\
&= \frac{1}{\sqrt{b_{H,k}}}  \E{ \left\Vert H_{i,k} -  \nabla_x^2 f(X_k)\right\Vert^2}^\sfrac{1}{2} \leq \frac{\sigma_H}{\sqrt{b_{H,k}}}. \nonumber 
\end{align}
Now injecting the last inequality and  \eqref{excpeHessgraderror} in \eqref{twotermsHik} yields that 
\[
\E{ \left\Vert \frac{1}{b_{H,k}} \sum_{i \in \calB_{H,k}} H_{i,k} -  \nabla_x^2 f(X_k) \right\Vert^3} 
\leq \left( 2 \sigma_H \sqrt{\frac{2e \log(n)}{b_{H,k}}}   
+ \frac{8e\log (n) \sigma_H}{b_{H,k}}  \right)^3.
\]
Imposing the left-hand side of the previous inequality to be less than
$\kappa_D \xi_k$ and using the concavity of the square root function then
yields that 
\begin{align*}
\frac{1}{\sqrt{b_{H,k}}} &\leq \frac{\sqrt{2e \log(n) + 8e \log(n)  \kappa_D^\sfrac{1}{3} \frac{\xi_k^\sfrac{1}{3}}{\sigma_H}} - \sqrt{2e \log(n)} }{8 e \log(n)} \\
&\leq \frac{8e\log(n) \kappa_D^\sfrac{1}{3} \xi_k^\sfrac{1}{3}}{12 e \sigma_H \log(n) \sqrt{2e \log(n)}} = \frac{2\kappa_D^\sfrac{1}{3} \xi_k^\sfrac{1}{3}}{3 \sigma_H \sqrt{2e \log(n)}}.
\end{align*}
Rearranging  the last inequality gives the bound \eqref{bHksamplebound}.
}

\section{Considered Datasets}\label{Datasets}
\begin{table}[!ht]
\centering
	\begin{tabular}{|l | r |r|}
		\hline
		Dataset & Samples & Features  \\ [0.5ex] 
		\hline
		a9a  & 32561 & 123  \\ [1.5ex] 
		ijcnn1  & 49990 & 22  \\ [1.5ex] 
		w8a  & 49749 & 300  \\ [1.5ex] 
		SUSY   & 5000000 & 18  \\ [1.5ex] 
		HIGGS  & 11000000 & 28  \\ 
		\hline
	\end{tabular}
\caption{\label{Datasettable}Datasets characterization, source:  LIBSVM\cite{LIBSVM} }
\end{table}

\section{Additional Results}\label{addfig}
\begin{figure}[!ht]
	\centering
	\includegraphics[width=7cm]{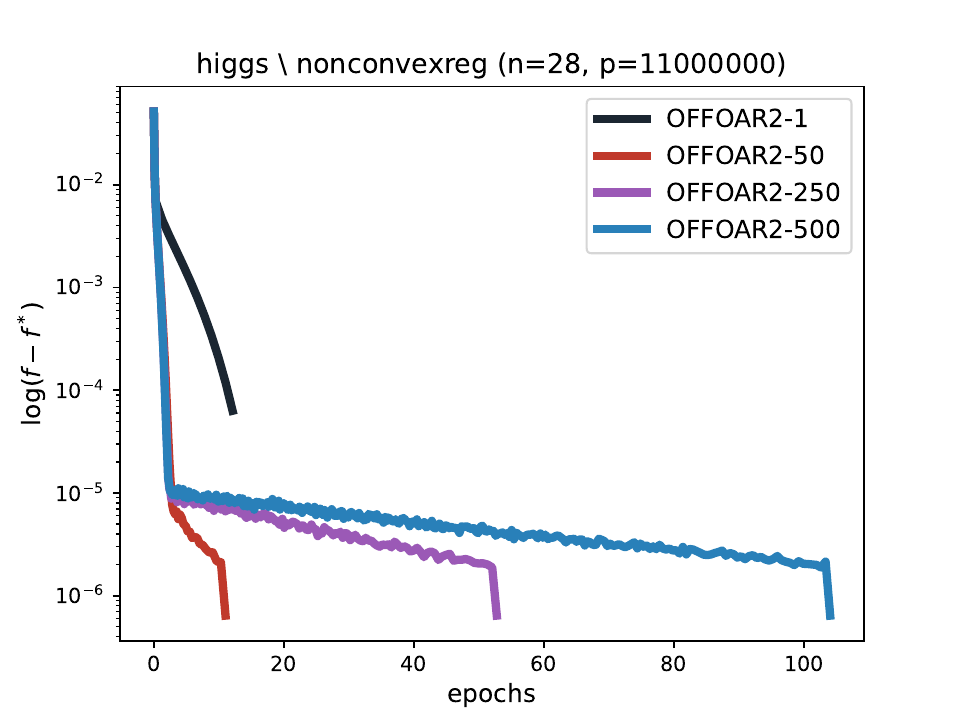}
	\hspace*{-5mm}
	\includegraphics[width=7cm]{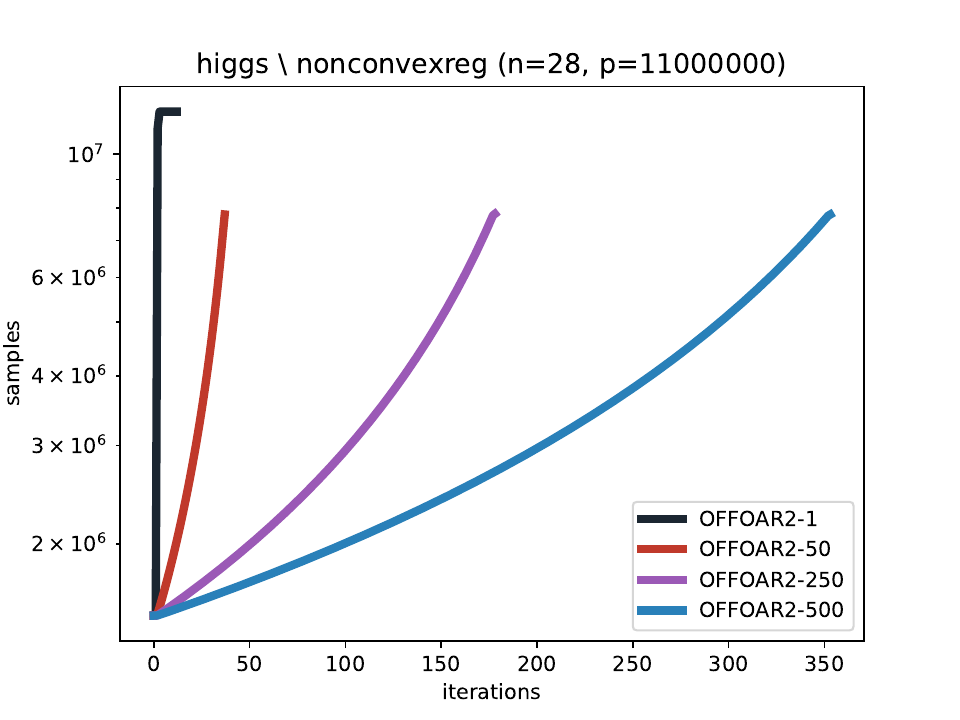}
	\hspace*{-5mm}
	\includegraphics[width=7cm]{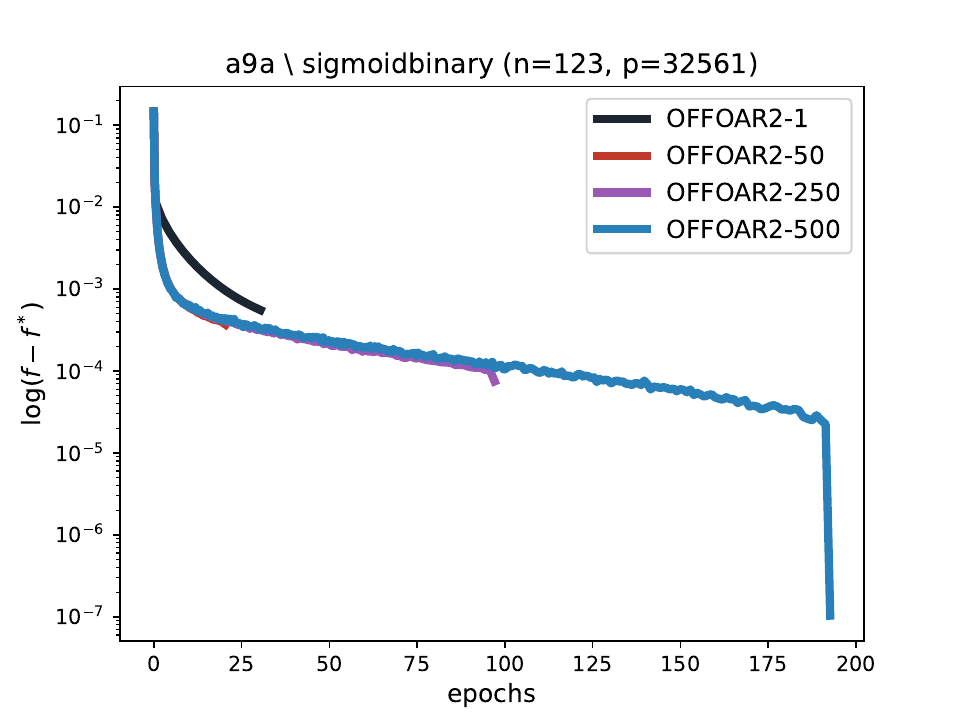}
	\hspace*{-5mm}
	\includegraphics[width=7cm]{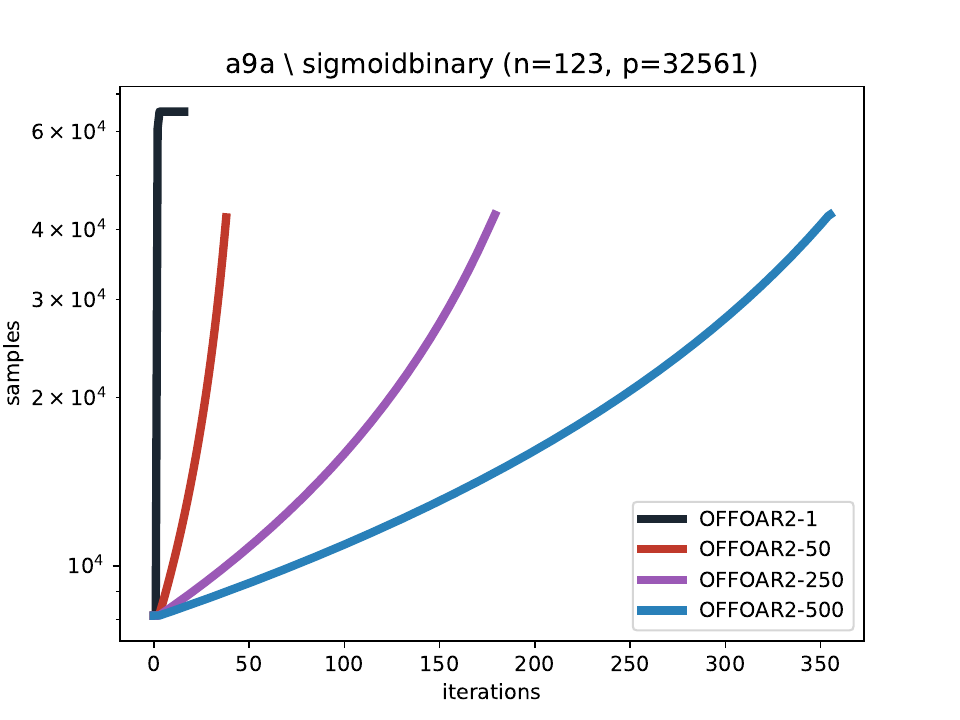}
	\hspace*{-5mm}
	\caption{\label{fig:objf2} Evolution of the loss function w.r.t the epochs and sampling behavior along iterations for specific problems}
\end{figure}

As shown in Figure~\ref{fig:objf2},  \texttt{OFFAR2-m} with large
$\texttt{m}$ may require a large number of epochs before achieving
convergence for some problems. From the sampling plots, we see that a
long-memory configuration may become an obstacle when the batch sizes
grow too slowly, thereby resulting in a substantial number of
iterations (and hence epochs) before achieving convergence.

\end{document}